\documentclass[11pt,twoside]{article}

\usepackage[T1]{fontenc}
\usepackage{amsmath,amsfonts,amssymb,amsthm,fullpage,bbm}
\usepackage{fancyhdr,graphicx,color}
\usepackage{epsfig}

\newtheorem{thm}{Theorem}
\newtheorem{lem}{Lemma}

\theoremstyle{definition}

\theoremstyle{remark}
\newtheorem{rem}{Remark}
\newenvironment{dem}{\textbf{Proof}~:\\}{\flushright$\blacksquare$\\}

\title{ UPPER LARGE DEVIATIONS FOR THE MAXIMAL FLOW IN FIRST PASSAGE PERCOLATION} 
\author{\Large Marie TH\'ERET}
\date{Laboratoire de Math\'ematiques, Universit\'e Paris Sud, B\^atiment
  425, 91405 Orsay, France}

\begin{document}
\maketitle

We consider the standard first passage percolation in $\mathbb{Z}^d$ for $d\geq2$ and we denote by $\smash{\phi_{n^{d-1},h(n)}}$ the maximal flow through the cylinder $\smash{]0,n]^{d-1} \times ]0,h(n)]}$ from its bottom to its top. Kesten proved a law of large numbers for the maximal flow in dimension three: under some assumptions, $\phi_{n^{d-1},h(n)}/n^{d-1}$ converges towards a constant $\nu$. We look now at the probability that $\smash{\phi_{n^{d-1},h(n)}/n^{d-1}}$ is greater than $\nu + \varepsilon$ for some $\varepsilon >0$, and we show under some assumptions that this probability decays exponentially fast with the volume $n^{d-1}h(n)$ of the cylinder. Moreover, we prove a large deviation principle for the sequence $\smash{(\phi_{n^{d-1},h(n)}/n^{d-1},n\in \mathbb{N})}$.


\section{Definitions and main results}

We will use generally the notations introduced in \cite{Kesten:StFlour} and \cite{Kesten:flows} but some changes will be done, for example to obtain independent objects. Let $d\geq 2$. We consider the graph $(\mathbb{Z}^{d}, \mathbb E ^{d})$ having for vertices $\mathbb Z ^{d}$ and for edges $\mathbb E ^{d}$ the set of all the pairs of nearest neighbors for the standard $L^{1}$ norm. With each edge $e$ in $\mathbb{E}^{d}$ we associate a random variable $t(e)$ with values in $\mathbb{R}^{+}$. We suppose that the family $(t(e), e \in \mathbb{E}^{d})$ is independent and identically distributed, with a common distribution function $F$. More formally, we take the product measure $\mathbb {P}$ on $\Omega= \prod_{e\in \mathbb{E}^{d}} [0, \infty[$, and we write its expectation $\mathbb{E}$. We interpret $t(e)$ as the capacity of the edge $e$; it means that $t(e)$ is the maximal amount of fluid that can go through the edge $e$ per unit of time. For a given realization $(t(e),e\in \mathbb{E}^{d})$ we denote by $\phi_{\vec{k},m} = \phi_{B}$ the maximal flow through the box
$$ B(\vec{k},m) \, = \, \prod _{i=1}^{d-1} ]0,k_{i}] \times ]0,m] \, ,$$
where $\vec{k}=(k_{1},...,k_{d-1}) \in \mathbb{Z}^{d-1}$, from its bottom
$$F_{0}\, = \, \prod _{i=1}^{d-1} ]0,k_{i}] \times \{0\} $$
to its top
$$F_{m}\, = \, \prod _{i=1}^{d-1} ]0,k_{i}] \times \{m\} \, .$$
Let us define this quantity properly. We recall that $\mathbb{E}^{d}$ is the set of the edges of the graph. An edge $e\in \mathbb{E}^{d}$ can be written $e=\langle x,y \rangle$, where $x$, $y\in \mathbb{Z}^{d}$ are the endpoints of $e$. The edges of $\mathbb{E}^{d}$ are unoriented, hence $\langle x,y \rangle = \langle y,x \rangle$. We will say that $e=\langle x,y\rangle$  is included in a subset $A$ of $\mathbb{R}^{d}$ ($e\subset A$) if the segment joining $x$ to $y$ (except possibly its extremities) is included in $A$. Now we define $\widetilde{\mathbb{E}}^{d}$ as the set of all the oriented edges, i.e., an element $\widetilde{e}$ in $\widetilde{\mathbb{E}}^{d}$ is an ordered pair of vertices. We denote an element $\widetilde{e} \in \widetilde{\mathbb{E}}^{d}$ by $\langle \langle x,y \rangle \rangle$, where $x$, $y \in \mathbb{Z}^{d}$ are the endpoints of $\widetilde{e}$ and the edge is oriented from $x$ towards $y$. We consider now the set $\mathcal{S}$ of all pairs of functions $(g,o)$, with $g:\mathbb{E}^{d} \rightarrow \mathbb{R}^{+}$ and $o:\mathbb{E}^{d} \rightarrow \widetilde{\mathbb{E}}^{d}$ such that $o(\langle x,y\rangle ) \in \{ \langle \langle x,y\rangle \rangle , \langle \langle y,x \rangle \rangle \}$, satisfying
\begin{itemize}
\item for each edge $e$ in $B$ we have
$$0 \,\leq\, g(e) \,\leq\, t(e) \,,$$
\item for each vertex $v$ in $B \smallsetminus F_{m}$ (remember that $F_{0}
  \cap B = \emptyset$) we have
$$ \sum_{e\in B\,:\, o(e)=\langle\langle v,\cdot \rangle \rangle} g(e) \,=\, \sum_{e\in B\,:\, o(e)=\langle\langle \cdot ,v \rangle \rangle} g(e) \,. $$ 
\end{itemize}
A couple $(g,o) \in \mathcal{S}$ is a possible stream in $B$: $g(e)$ is the amount of fluid that goes through the edge $e$, and $o(e)$ gives the direction in which the fluid goes through $e$. The first condition on $(g,o)$ expresses only the fact that the amount of fluid that can go through an edge is bounded by its capacity. The second one is a balance equation: it means that there is no loss of fluid in the cylinder. With each possible stream we associate the corresponding flow
$$ flow (g,o) \,=\, \sum_{ u \in B\smallsetminus F_{m} \,,\,  v \in F_{m} \,:\, \langle u,v\rangle \in \mathbb{E}^{d}} g(\langle u,v\rangle) \mathbb{I}_{o(\langle u,v\rangle) = \langle\langle u,v \rangle\rangle} - g(\langle u,v\rangle) \mathbb{I}_{o(\langle u,v\rangle) = \langle\langle v,u \rangle\rangle} \,. $$
This is the amount of fluid that crosses the cylinder $B$ if the fluid respects the stream $(g,o)$. The maximal flow through the cylinder $B$ from its bottom to its top is the supremum of this quantity over all possible choices of streams in $\mathcal{S}$
$$ \phi_{B} \,=\, \phi_{\vec{k},m} \,=\, \sup \, \{ flow (g,o) \,:\, (g,o) \in \mathcal{S} \} \,.$$
If $\phi_{B} = flow (g,o)$ we say that the stream $(g,o)$ realizes the flow $\phi_{B}$.

Kesten proved in 1987 the following law of large numbers for the maximal flow in dimension $3$ (see \cite{Kesten:flows}):
\begin{thm}
\label{lgn}
We consider a cylinder $B((k,l),m)$ such that $\lim_{k \geq l \rightarrow \infty} m(k,l) =\infty$ in such a way that for some $\delta >0$ we have
$$ \lim_{k\geq l \rightarrow \infty} \frac{\ln m(k,l)}{k^{1-\delta}} \,=\, 0 \,. $$
There exists a positive $p_{0}$ with the following property: If $F$
satisfies $F(0) < p_{0}$ and $\int_{[0,\infty[} e^{\theta x} dF(x)$ is finite for some positive $\theta$, then there exists a constant $\nu (F) < \infty$ such that
$$ \lim_{k,l \rightarrow \infty} \frac{\phi_{(k,l),m}}{kl} \,=\, \nu  \qquad with\,\, probability\,\, one\,\, and\,\, in\,\, L^{1} \,. $$
\end{thm}
Actually, the constant $\nu(F)$ is defined as the limit of another object
under weaker assumptions on $F$ (see \cite{Kesten:flows} and (\ref{defnu}) in the
next section), and we rely on this definition to state the following result. We are now interested in the deviations of the rescaled flow from its typical behavior. We will show two results in dimensions $d\geq 2$. The first one states the existence of a limit, and some of its properties.

\begin{thm}
\label{limite}
We consider the maximal flow $\phi_{(n,...,n),h(n)}$ through the cylinder $B((n,...,n),h(n))$, where the function $h: \mathbb N \rightarrow \mathbb N$ satisfies
$$ \lim_{n\rightarrow \infty} \frac{h(n)}{\ln{n}} \, = \, \infty  \, .$$
For every $\lambda$ in $\mathbb R ^{+}$, the limit
$$ \psi(\lambda) \, = \, \lim_{n\rightarrow \infty} - \frac{1}{n^{d-1}h(n)} \ln{ \mathbb P \left[\phi_{(n,...,n),h(n)} \geq \lambda n^{d-1}\right] }  $$
exists and is independent of $h$.
Moreover $\psi$ is convex on $\mathbb R ^{+}$, finite and continuous on the
set $\{ \, \lambda \, | \, F([\lambda , +\infty [)  >0 \, \}$. If $
    \int_{[0,+\infty[} x dF(x) $ is finite, then
    $\psi$ vanishes on $[0,\nu]$, where $\nu$ is defined in (\ref{defnu}). If $ \int_{[0,+\infty [} e^{\theta x}
        dF(x)$ is finite for some positive $\theta$, then $\psi$ is positive on $]\nu,+\infty[$.
\end{thm}

We say that a sequence $(X_{n}, n\in \mathbb{N})$ of random variables with
values in $D \subset \mathbb{R}$ satisfies a large deviation principle
with speed $v(n)$ and governed by the rate function $\mathcal{I}$ if and
only if
\begin{itemize}
\item for any closed subset $\mathcal{F} \subset D$, we have
$$ \limsup _{n\rightarrow \infty} \frac{1}{v(n)} \ln \mathbb{P} \left[ X_{n} \in \mathcal{F} \right] \,\leq \, -\inf_{\mathcal{F}} \mathcal{I} \,, $$
\item for any open subset $\mathcal{O} \subset D$, we have
$$\liminf _{n\rightarrow \infty} \frac{1}{v(n)} \ln \mathbb{P} \left[ X_{n} \in \mathcal{O} \right] \,\geq \, -\inf_{\mathcal{O}} \mathcal{I} \,.  $$
\end{itemize}
Now, with the help of the function $\psi$, we can state the following large
deviation principle for the rescaled flow:
\begin{thm}
\label{pgd}
Let $h: \mathbb N \rightarrow \mathbb N$ be such that
$$\lim_{n\rightarrow \infty}\frac{h(n)}{\ln n} = \infty \,.$$
If there exists a positive $\theta$ such that
$$\int_{[0, +\infty[}  e^{\theta x} dF(x) \, < \, \infty \, ,$$
then the sequence
$$\left(\frac{\phi_{(n,...,n),h(n)}}{n^{d-1}} \right)_{n \in \mathbb N}$$
satisfies a large deviation principle, with speed $n^{d-1}h(n)$, and governed by the good rate function $\psi$.
\end{thm}

\begin{rem}
When the capacity $t$ of an edge is bounded, we do not need the condition
$$ \lim_{n\rightarrow \infty} \frac{h(n)}{\ln n} \,=\, +\infty \,,$$
the results hold under the weaker condition
$$ \lim_{n\rightarrow \infty}h(n) \,=\, +\infty \,. $$
Actually, the role of the condition $\lim_{n \rightarrow \infty} h(n)/\ln n = +\infty$ is not fully understood yet. For example, when  $t$ is equal in law to the absolute value of a Gaussian variable this condition can also be replaced by $\lim_{n\rightarrow +\infty} h(n) = +\infty$. Unfortunately we could not find satisfying sufficient conditions (in particular on the moments of the law of $t$) to get rid of the condition  $\lim_{n \rightarrow \infty} h(n)/\ln n = +\infty $.
\end{rem}

A special aspect of the proof of theorem \ref{limite} is the use of a discrete
version of the model. Indeed, we are confronted with a combinatorial
problem: we need to look at boundary conditions for streams to glue
together streams in different cylinders, but when the capacity of an edge
takes its values in $\mathbb{R}^{+}$ we cannot count the number of possible
boundary conditions. Our strategy is to consider a discrete approximation
of the capacity of the edges and the corresponding maximal flow. We work
with these objects. To handle the boundary conditions, we use a technique
introduced by Chow and Zhang \cite{Chow-Zhang}. We finally compare the real maximal flow to this approximation.


\section{Max-flow min-cut theorem}
\label{*}

It is difficult to work with the expression of the maximal flow that we have seen in the previous part, this is the reason why we will use the max-flow min-cut theorem to express the maximal flow differently. First we need some definitions. A path on a graph ($\mathbb{Z}^{d}$ for example) from $v_{0}$ to $v_{n}$ is a sequence $(v_{0}, e_{1}, v_{1},..., e_{n}, v_{n})$ of vertices $v_{0},..., v_{n}$ alternating with edges $e_{1},..., e_{n}$ such that $v_{i-1}$ and $v_{i}$ are neighbors in the graph, joined by the edge $e_{i}$, for $i$ in $\{1,..., n\}$. Two paths are said disjoint if they have no common edge. A set $E$ of edges of $B(\vec{k},m)$ is said to separate $F_{0}$ from $F_{m}$ in $B(\vec{k},m)$ if there is no path from $F_{0}$ to $F_{m}$ in $B(\vec{k},m) \smallsetminus E$. We call $E$ an $(F_{0},F_{m})$-cut if $E$ separates $F_{0}$ from $F_{m}$ in $B(\vec{k},m)$ and if no proper subset of $E$ does. With each set $E$ of edges we associate the variable
$$ V(E)\, = \, \sum_{e\in E} t(e) \, .$$
The max-flow min-cut theorem (see \cite{Bollobas}) states that
$$ \phi_{B} \, = \, \min \{ \, V(E) \, | \, E \,\, is \,\, an \,\, (F_{0},F_{m})-cut \,\, in \,\, B \, \} \, .$$

\begin{rem}
\label{menger}
In the special case where $t(e)$ belongs to $\{0,1\}$, let us consider the graph obtained from the initial graph (not necessarily $\mathbb{Z}^{d}$) by removing all the edges $e$ with $t(e)=0$. Menger's theorem (see \cite{Bollobas}) states that the minimal number of edges in $B(\vec{k},m)$ that have to be removed from this graph to disconnect $F_{0}$ from $F_{m}$ is exactly the maximal number of disjoint paths that connect $F_{0}$ to $F_{m}$. By the max-flow min-cut theorem, it follows immediately that the maximal flow in the initial graph through $B$ from $F_{0}$ to $F_{m}$ is exactly the maximal number of disjoint open paths from $F_{0}$ to $F_{m}$, where a path is open if and only if the capacity of all its edges is one. Such a set of $\phi_{B}$ disjoint open paths from $F_{0}$ to $F_{m}$ corresponds obviously to a stream (g,o):
\begin{itemize}
\item$ g(e) \, = \, \left\{ \begin{array}{ll} 1 &\qquad if\,\,e\,\,belongs\,\,to\,\,one\,\, of\,\,these\,\,paths \\ 0 & \qquad otherwise\,, \end{array} \right.$
\item $o(e) \,=\, \left\{ \begin{array}{ll} \langle\langle x,y\rangle\rangle &  \qquad if \,\,e=\langle x,y\rangle \,\, is \,\,crossed\,\,from\,\,x\,\,to\,\,y\,\,by\,\,one\,\,of\,\,these\,\,paths \\ 
\langle\langle y,x\rangle\rangle &  \qquad if \,\,e=\langle x,y\rangle \,\, is \,\,crossed\,\,from\,\,y\,\,to\,\,x\,\,by\,\,one\,\,of\,\,these\,\,paths \\ 
\hat{o} (e) & \qquad otherwise\,,  \end{array} \right.$
\end{itemize}
where $\hat{o}$ is some determined orientation ($\hat{o}(\langle x,y\rangle ) \in \{ \langle \langle x,y\rangle\rangle , \langle\langle y,x\rangle\rangle \}$) which does not matter. The stream $(g,o)$ realizes the maximal flow $\phi_{B}$ (whatever $\hat{o}$).
\end{rem}

We come back to the general case. We will also need the definition of a cut over a hyper-rectangle. Let $S= \prod_{i=1}^{d-1} ]k_{i},l_{i}]$ be a hyper-rectangle, with $k_{i} \leq l_{i}$, $k_{i}$, $l_{i}$ in $\mathbb Z$. We say that a set $E$ of edges in $S \times \mathbb R$ separates $-\infty$ from $+\infty$ over $S$ if there exists no path in $(S\times \mathbb R) \smallsetminus E$ from $S \times \{-N\}$ to $S\times \{+N\}$ for some $N>0$. Similarly, we call $E$ a cut over $S$ if $E$ separates $-\infty$ from $+\infty$ over $S$, but no proper subset of $E$ does. Let $\partial ^{in} (S\times \mathbb R)$ be the inner vertex boundary of the cylinder $S\times \mathbb R$
$$ \partial ^{in} (S\times \mathbb R) \,=\, \{\, x\in S\times \mathbb R \,|\, \exists y\notin S\times \mathbb R \,,\,\, \langle x,y \rangle \in \mathbb{E}^{d} \, \} \,. $$
We define the corresponding set of edges
$$ \mathbb{E}(\partial ^{in}( S\times \mathbb R)) \,=\, \{\langle
x,y\rangle \,|\, x,y\in \partial ^{in} ( S\times \mathbb R)  \}\,.$$
We say that an edge $e$ is vertical if $e=\langle x,
x+(0,...,0,1) \rangle$; $e$ is said horizontal otherwise. We denote by $(*)$ the condition on $E$
$$ (*) \qquad E \cap \mathbb{E}( \partial ^{in}( S\times \mathbb R) ) \, \subset \, \{ e\in \mathbb{E}^{d} \,|\, e\,\, is \,\, vertical\,,\,\, e\subset \mathbb{R}^{d-1}\times [0,1] \} \,,  $$
which means in a way to say that the boundary of $E$ is fixed on the perimeter of the rectangle $S\times \{0\}$. We define the variable $\tau$ by
$$ \tau(S) \, = \, \inf \, \{ \, V(E) \, |\, E \,\, is \,\, a \,\, cut \,\, over \,\, S \,\, and \,\, E \,\, satisfies \,\, (*) \, \} \, .$$
For simplicity, we denote by $\tau_{k^{d-1}}$ the variable $\tau(
          ]0,k]^{d-1})$. If $S_{1}$, $S_{2}$ are two disjoint
        hyper-rectangles having a common side (so $S_{1} \cup S_{2}$ is an hyper-rectangle too), then we
        have
$$ \tau(S_{1} \cup S_{2}) \, \leq \,  \tau(S_{1}) + \tau (S_{2}) \,.  $$
Indeed if $E_{1}$ (respectively $E_{2}$) is a cut over $S_{1}$
(respectively $S_{2}$) satisfying $(*)$ for $S_{1}$ (respectively $S_{2}$)
then $E_{1}$ and $E_{2}$ are both pinned at the boundary between
$S_{1} \times \{ 0 \}$ and $S_{2} \times \{ 0 \}$ because they both satisfy
$(*)$, so they can be glued together and $E_{1} \cup E_{2}$ separates $-\infty$ from $+\infty$ over $S_{1}
\cup S_{2}$. By a subadditive argument (see \cite{Ackoglu}), the following limit exists almost surely
\begin{equation}
\label{defnu}
\nu(F) \, = \, \lim_{k\rightarrow \infty} \frac{\tau_{k^{d-1}}}{k^{d-1}} \,
,
\end{equation}
where we know that $\nu(F)$ is a constant almost surely thanks to
Kolmogorov's $0-1$ law. We will denote it by $\nu$ when no doubt about $F$ is possible. This is the ``$\nu$'' in theorems \ref{lgn} and \ref{limite}.


\section{Proof of Theorem \ref{limite}}

We take $h: \mathbb N \rightarrow \mathbb N$ such that
$$\lim_{n\rightarrow \infty} h(n)\,=\, +\infty \,.$$
We will see during the proof where we need the stronger condition
$$\lim_{n\rightarrow \infty} \frac{h(n)}{\ln n}\,=\, +\infty \,.$$
We will need to describe how the fluid goes in and out of a cylinder in
order to glue together two cylinders without loosing any flow. The problem
is that we need too much information to describe this precisely. The
feature of the proof is to consider a discrete approximation of the
capacity of the edges (see section \ref{discretisation}), to work with this
discrete model (sections \ref{limitediscrete}, \ref{convexite} and
\ref{continuite}) and then to compare it to the original one (section \ref{limitegen}). The method used to prove the existence of the limit was developed in \cite{Chow-Zhang}. We study then the properties of $\psi$ as in \cite{Cerf:StFlour}. 

\subsection{Discrete version}
\label{discretisation}

Let $k \in \mathbb N$ (we will choose it later). We associate with $(t(e), e \in \mathbb Z ^{d})$ a new family of independent and identically distributed variables $(t^{k}(e), e \in \mathbb Z ^{d})$ by setting
$$ \forall e\in \mathbb{E}^{d} \qquad t^{k}(e) \, = \, \lfloor kt(e)  \rfloor \times \frac{1}{k} \, ,$$
and we denote by $\phi^{k}$ the maximal flow corresponding to these new variables.

\begin{figure}[ht!]
\centering
\begin{picture}(0,0)%
\epsfig{file=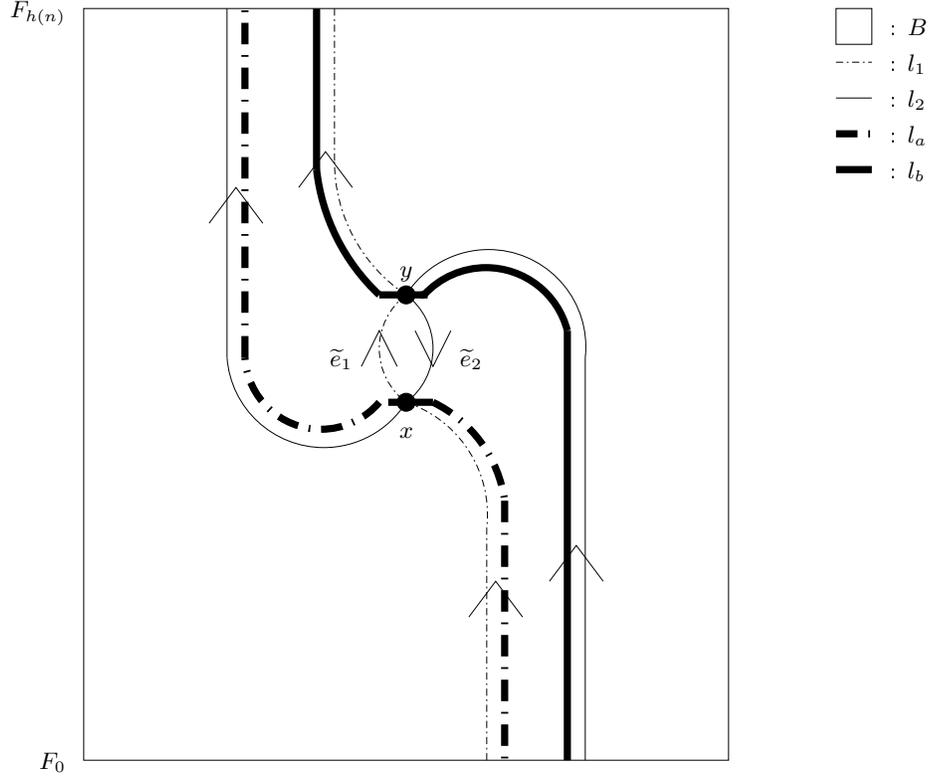}%
\end{picture}%
\setlength{\unitlength}{2960sp}%
\begingroup\makeatletter\ifx\SetFigFont\undefined%
\gdef\SetFigFont#1#2#3#4#5{%
  \reset@font\fontsize{#1}{#2pt}%
  \fontfamily{#3}\fontseries{#4}\fontshape{#5}%
  \selectfont}%
\fi\endgroup%
\begin{picture}(8500,6514)(1343,-6494)
\put(5101,-2311){\makebox(0,0)[b]{\smash{{\SetFigFont{9}{10.8}{\rmdefault}{\mddefault}{\updefault}{\color[rgb]{0,0,0}$y$}%
}}}}
\put(5551,-3061){\makebox(0,0)[lb]{\smash{{\SetFigFont{9}{10.8}{\rmdefault}{\mddefault}{\updefault}{\color[rgb]{0,0,0}$\widetilde{e}_{2}$}%
}}}}
\put(4651,-3061){\makebox(0,0)[rb]{\smash{{\SetFigFont{9}{10.8}{\rmdefault}{\mddefault}{\updefault}{\color[rgb]{0,0,0}$\widetilde{e}_{1}$}%
}}}}
\put(9151,-286){\makebox(0,0)[lb]{\smash{{\SetFigFont{9}{10.8}{\rmdefault}{\mddefault}{\updefault}{\color[rgb]{0,0,0}: $B$}%
}}}}
\put(9151,-586){\makebox(0,0)[lb]{\smash{{\SetFigFont{9}{10.8}{\rmdefault}{\mddefault}{\updefault}{\color[rgb]{0,0,0}: $l_{1}$}%
}}}}
\put(9151,-886){\makebox(0,0)[lb]{\smash{{\SetFigFont{9}{10.8}{\rmdefault}{\mddefault}{\updefault}{\color[rgb]{0,0,0}: $l_{2}$}%
}}}}
\put(9151,-1186){\makebox(0,0)[lb]{\smash{{\SetFigFont{9}{10.8}{\rmdefault}{\mddefault}{\updefault}{\color[rgb]{0,0,0}: $l_{a}$}%
}}}}
\put(9151,-1486){\makebox(0,0)[lb]{\smash{{\SetFigFont{9}{10.8}{\rmdefault}{\mddefault}{\updefault}{\color[rgb]{0,0,0}: $l_{b}$}%
}}}}
\put(2251,-6436){\makebox(0,0)[rb]{\smash{{\SetFigFont{9}{10.8}{\rmdefault}{\mddefault}{\updefault}{\color[rgb]{0,0,0}$F_{0}$}%
}}}}
\put(2251,-136){\makebox(0,0)[rb]{\smash{{\SetFigFont{9}{10.8}{\rmdefault}{\mddefault}{\updefault}{\color[rgb]{0,0,0}$F_{h(n)}$}%
}}}}
\put(5101,-3661){\makebox(0,0)[b]{\smash{{\SetFigFont{9}{10.8}{\rmdefault}{\mddefault}{\updefault}{\color[rgb]{0,0,0}$x$}%
}}}}
\end{picture}%
\caption{Description of paths in $G$}
\label{courant}
\end{figure}

Let us consider for a brief moment the graph $G$ obtained by replacing each edge $e$ by $p$ edges $\widetilde{e}_{1},...,\widetilde{e}_{p}$, where $p=k t^{k}(e)$. In this new graph the capacity of each edge is simply one. The remark \ref{menger} also holds for $G$: the maximal flow $\phi^{G}_{B}$ for $G$ from $F_{0}$ to $F_{h(n)}$ in $B=B((n,...,n),h(n))$ is exactly the maximal number of disjoint paths connecting $F_{0}$ to $F_{h(n)}$ in $G$. We have seen that we can associate with each such family of $\phi_{B}$ disjoint paths in $B$ a stream $(\widetilde{g},\widetilde{o})$ in $G$ that realizes $\phi_{B}^{G}$. Actually, we can always reduce to the case where $\widetilde{o}(\widetilde{e}_{1}) = \widetilde{o}(\widetilde{e}_{2})$ if the edges $\widetilde{e}_{1}$ and $\widetilde{e}_{2}$ are replacing in $G$ the same edge $\langle x,y\rangle \in \mathbb{E}^{d}$. Indeed, if for such edges $\widetilde{e}_{1}$ and $\widetilde{e}_{2}$ we have $\widetilde{g}(\widetilde{e}_{1})=\widetilde{g}(\widetilde{e}_{2})=1$ and $\widetilde{o}(\widetilde{e}_{1}) \neq \widetilde{o}(\widetilde{e}_{2})$, we know that there exists a path $l_{1}$ (respectively $l_{2}$) from $F_{0}$ to $F_{h(n)}$ going through $\widetilde{e}_{1}$ (respectively $\widetilde{e}_{2}$) and crossing this edge from $x$ to $y$ (respectively from $y$ to $x$). We can create two new disjoint paths in $G$, $l_{a}$ which is equal to $l_{1}$ from $F_{0}$ to $x$ and to $l_{2}$ from $x$ to $F_{h(n)}$, and $l_{b}$ which is equal to $l_{2}$ from $F_{0}$ to $y$ and to $l_{1}$ from $y$ to $F_{h(n)}$, that can replace $l_{1}$ and $l_{2}$ in the set of $\phi_{B}$ disjoint paths (see figure \ref{courant}). The corresponding stream $(\widetilde{g}',\widetilde{o}')$ is equal to $(\widetilde{g},\widetilde{o})$ except in $\widetilde{e}_{1}$ and $\widetilde{e}_{2}$ where we have $\widetilde{g}'(\widetilde{e}_{1}) = \widetilde{g}'(\widetilde{e}_{2}) = 0$ and $\widetilde{o}'(\widetilde{e}_{1}) = \widetilde{o}'(\widetilde{e}_{2}) = \hat{o}(\langle x,y\rangle )$.
Then we just have to deal with the case
$\widetilde{g}(\widetilde{e}_{1})=1$ and
$\widetilde{g}(\widetilde{e}_{2})=0$. The definition of $\hat{o}$ is
arbitrary, we can change the orientation of $\hat{o} (\widetilde{e}_{2})$
to have $\widetilde{o} (\widetilde{e}_{2}) = \hat{o} (\widetilde{e}_{2}) =
\widetilde{o} (\widetilde{e}_{1})$. We obtain thus a stream such that $\widetilde{o}(\widetilde{e}_{1}) = \widetilde{o}(\widetilde{e}_{2})$ if the edges $\widetilde{e}_{1}$ and $\widetilde{e}_{2}$ are replacing in $G$ the same edge $\langle x,y\rangle \in \mathbb{E}^{d}$ . Moreover we can assume that each path of a family of $\phi_{B}$ disjoint open paths has only its first vertex in $F_{0}$ and its last vertex in $F_{h(n)}$, otherwise we can restrict the path to obtain such a path. Thanks to a good choice of $\hat{o}$, we can thus suppose that if $\widetilde{e} \subset B$ has one endpoint $x$ in $F_{h(n)}$ (respectively $F_{0}$) and one endpoint $y$ not in $F_{h(n)}$ (respectively $F_{0}$) then $\widetilde{o}(\widetilde{e})= \langle \langle y,x \rangle\rangle$ (respectively $\widetilde{o}(\widetilde{e})= \langle \langle x,y \rangle\rangle$), and if $\widetilde{e}$ has both endpoints in $F_{h(n)}$ then $\widetilde{g}(\widetilde{e})=0$.

Coming back to the graph $\mathbb{Z}^{d}$, we remark that the maximal flow $\phi_{B}^{k}$ between $F_{0}$ and $F_{h(n)}$ in $B$ is equal to $\phi_{B}^{G}/k$, and it can be realized by the stream $(g,o)$ defined as follows. Let $e$ be an edge of $\mathbb{E}^{d}$. If there is no edge in $G$ associated with $e$, we set $g(e)=0$ and $o(e)=\hat{o}(e)$. Otherwise, we define
$$ g(e) \,=\, \sum_{\widetilde{e} \sim e} \frac{\widetilde{g}(\widetilde{e})}{k}$$
where the sum is over the edges $\widetilde{e}$ that replace $e$ in $G$,
and $o(e) = \widetilde{o}(\widetilde{e})$ for some edge $\widetilde{e}$
associated with $e$ (recall that if $\widetilde{e}_{1} \sim e $ and
$\widetilde{e}_{2} \sim e$ then
$\widetilde{o}(\widetilde{e}_{1})=\widetilde{o}(\widetilde{e}_{2})$). We
will call such a stream, built from the graph $G$, a discrete stream. A
discrete stream has three particular properties: $g$ takes its values in
$\mathbb{N}/k$, $o(\langle x,y\rangle)=\langle\langle x,y \rangle\rangle$
as soon as we have $x\in F_{0}$ and $y\in B$ ($y \notin F_{0}$) or $y\in
F_{h(n)}$ and $x\in B\smallsetminus F_{h(n)}$, and $g(e)=0$ if $e$ has both endpoints in $F_{h(n)}$.

Let $\lambda$ be in $\mathbb{R}^{+}$. For a discrete stream $(g,o)$ and
$h\in \mathbb{Z}$ we define the truncated projection of $g$ on the vertical
edges that intersect the hyper-plane $\{(x_{1},...,x_{d}) \in
\mathbb{R}^{d} \,|\, x_{d} = h+1/2 \}$ by
$$ \forall x \in \mathbb{Z}^{d-1} \cap\, ]0,n]^{d-1} \qquad \pi_{h}^{\lambda,n}(g,x) \,=\, g\left( \langle (x,h),(x,h+1) \rangle \right) \wedge \left( \lfloor \lambda n^{d-1} \rfloor +1 \right) \,.$$
Thanks to the properties of discrete streams, we can state the following lemma:
\begin{lem}[Junction of two boxes]
\label{junction}
Let $B_{1}=]0,n]^{d-1} \times ]0,h(n)]$, $B_{2}=]0,n]^{d-1} \times
]h(n),2h(n)]$. If there exist a discrete stream $(g_{1},o_{1})$ in $B_{1}$
and a discrete stream $(g_{2},o_{2})$ in $B_{2}$ such that
$$ flow(g_{1},o_{1}) \,\geq\, \lambda n^{d-1} \qquad and \qquad  flow(g_{2},o_{2}) \,\geq\, \lambda n^{d-1} $$
and
$$\forall x\in \mathbb{Z}^{d-1} \cap\, ]0,n]^{d-1} \qquad  \pi_{h(n)-1}^{\lambda,n}(g_{1},x) \,=\,
    \pi_{h(n)}^{\lambda,n}(g_{2},x) \,, $$
then $\phi_{B_{1}\cup B_{2}}^{k} \geq \lambda n^{d-1}$.
\end{lem}

\begin{dem}
To prove this lemma, we consider two cases:
\begin{enumerate}
\item if for all $e=\langle (x,h(n)-1),(x,h(n)) \rangle$ with $x\in
  \mathbb{Z}^{d-1}\cap\, ]0,n]^{d-1}$ we have $g_{1}(e) \leq \lambda n^{d-1}$, then we can define the following discrete stream $(g_{tot},o_{tot})$:
\begin{itemize}
\item $g_{tot}(e) \,=\, \left\{ \begin{array}{ll} g_{1}(e) & \qquad if\,\,e\subset B_{1} \\
g_{2}(e) & \qquad if\,\,e\subset B_{2} \\
0 & \qquad otherwise \,, \end{array} \right.$
\item $o_{tot}(e) \,=\, \left\{ \begin{array}{ll} o_{1}(e) & \qquad if\,\,e\subset B_{1} \\
o_{2}(e) & \qquad if\,\,e\subset B_{2} \\
\hat{o}(e) & \qquad otherwise \,, \end{array} \right. $
\end{itemize}
where $\hat{o}$ is still some arbitrarily determined orientation. We can check that $(g_{tot},o_{tot})$ is a discrete stream thanks to the properties of $(g_{1},o_{1})$ and $(g_{2},o_{2})$, in particular if $e_{1}=\langle (x,h(n)-1),(x,h(n)) \rangle$ and $e_{2}=\langle (x,h(n)),(x,h(n)+1)\rangle$ with $x\in ]0,n]^{d-1}$, we have $g(e_{1})=g(e_{2})$, $g(e)=0$ for all others edges $e=\langle (x,h(n)), \cdot \rangle$, $o(e_{1})=\langle\langle(x,h(n)-1),(x,h(n)) \rangle\rangle$ and $o(e_{2})=\langle\langle (x,h(n)),(x,h(n)+1) \rangle\rangle$, hence the balance equation is satisfied. Moreover $flow(g_{tot},o_{tot}) = flow(g_{1},o_{1})  = flow(g_{2},o_{2}) $ so $\phi_{B_{1}\cup B_{2}}^{k} \geq \lambda n^{d-1}$.
\item Suppose there exists an edge $e = \langle (x,h(n)-1),(x,h(n))  \rangle$ such that $g_{1}(e)>\lambda n^{d-1}$. The discrete stream $(g_{1},o_{1})$ corresponds to $k \times flow(g_{1},o_{1})$ disjoint paths from $F_{0}$ to $F_{h(n)}$ in $B_{1}$ for the modified graph $G$. The inequality $g_{1}(e)>\lambda n^{d-1}$ implies that at least $q= \lceil \lambda n^{d-1} k \rceil $ of these paths, that we will denote by $l_{1},...,l_{q}$, go out of $B_{1}$ through $e$. The equality $\pi_{h(n)-1}^{\lambda,n}(g_{1},x) = \pi_{h(n)}^{\lambda,n}(g_{2},x)$ implies that $g_{2}(f) >\lambda n^{d-1}$ where $f=\langle (x,h(n)),(x,h(n)+1) \rangle$. By the same argument, we can find at least $q$ disjoint paths $l_{1}',...,l_{q}'$ from $F_{h(n)}$ to $F_{2h(n)}$ in $B_{2}$ for $G$, all going in $B_{2}$ through the edge $f$. Now we can glue together these $q$ paths $l_{1},...,l_{q}$ in $B_{1}$ with the $q$ paths $l_{1}',...,l_{q}'$ in $B_{2}$ because $e$ and $f$ are adjacent. This way we obtain $q$ disjoint paths from $F_{0}$ to $F_{2h(n)}$ in $B_{1} \cup B_{2}$ for $G$, and by considering the corresponding discrete flow $(g_{tot},o_{tot})$ in the initial graph we obtain $\phi_{B_{1}\cup B_{2}}^{k} \geq flow (g_{tot},o_{tot}) = \lceil \lambda n^{d-1} k \rceil /k$.
\end{enumerate}
\end{dem}

We define the boundary conditions of the discrete stream $(g,o)$ in the cylinder $B((n,...,n),h(n))$ as
\begin{align*}
\Pi ^{\lambda,n} (g) & \,=\, \left( \Pi^{\lambda,n}_{1}(g),\Pi^{\lambda,n}_{2}(g)  \right) \\
& \,=\, \left( \left(\pi^{\lambda,n}_{0}(g,x), x\in \mathbb{Z}^{d-1}\cap\,
]0,n]^{d-1}\right), \left(\pi^{\lambda,n}_{h(n)-1}(g,x), x\in
    \mathbb{Z}^{d-1} \cap\, ]0,n]^{d-1}\right) \right) \,.\\
\end{align*}
The number $N_{\lambda,n}^{k}$ of possible boundary conditions for discrete streams satisfies
$$ N_{\lambda,n}^{k} \, \leq \, \left( k \left( \lfloor \lambda n^{d-1} \rfloor +1 \right) +1 \right)^{2n^{d-1}} \, .$$


\subsection{Existence of the limit for $\phi^{k_{n}}_{n^{d-1},h(n)}$}
\label{limitediscrete}

In this section, we will prove the existence of the limit appearing in the theorem \ref{limite} with $\phi^{k}$ instead of $\phi$. We denote $\phi_{(n,...,n),h(n)}$ by $\phi_{n^{d-1},h(n)}$, and we define
$$ \mu \, = \, \sup \{ \, \lambda \, | \, F([0,\lambda[)<1 \, \} \, .$$
We will prove the following result:
\begin{thm}
\label{thmdiscret}
For every pair $(h,(\chi_{n}, n\in \mathbb{N}) ) $ with $h :
\mathbb{N} \rightarrow \mathbb{N}$ a function such that $\lim_{n\rightarrow + \infty}
h(n) = + \infty$ and $(\chi_{n}, n\in \mathbb{N})$ a non-decreasing sequence of integers such that $\lim_{n\rightarrow +\infty}\chi_{n} = +\infty
$, satisfying
\begin{equation}
\label{h-chi}
\lim_{n\rightarrow +\infty} \frac{\chi_{n} \ln n}{h(n)} \, = \, 0 \,,
\end{equation}
for every $\lambda$ in $\mathbb{R}^{+} \smallsetminus \{ \mu \}$ (or in
$\mathbb{R}^{+} $ if $\mu$ is infinite), the limit
$$ \widetilde{\psi} (\lambda, h, (\chi_{n})) \, = \, \lim_{n\rightarrow
  +\infty} - \frac{1}{n^{d-1}h(n)} \ln \mathbb{P} \left[
  \phi_{n^{d-1},h(n)}^{2^{\chi_{n}}} \geq \lambda n^{d-1}  \right]  $$
exists. Moreover $\widetilde{\psi}$ is independent of such a pair
  $(h,(\chi_{n}))$, i.e. if $(h,(\chi_{n}))$ and
  $(h',(\chi'_{n}))$ satisfy all the previous conditions, then
  $\widetilde{\psi} (\lambda, h, (\chi_{n})) = \widetilde{\psi} (\lambda, h',
  (\chi'_{n}))$ for all $\lambda$ in $\mathbb{R}^{+} \smallsetminus \{ \mu \}$
  (or in $\mathbb{R}^{+}$). We will thus denote this limit by $\widetilde{\psi} (\lambda)$.
\end{thm}
We now prove theorem \ref{thmdiscret} by considering different cases.

$\bullet \,\, \lambda>\mu$ : Then
$$ \forall k \in \mathbb{N}\,,\,\, \forall n \in \mathbb{N} \qquad \mathbb P \left[\phi^{k}_{n^{d-1},h(n)} \geq \lambda n^{d-1} \right] \, = \, 0 \,,$$
so for every sequence $(\chi_{n})$ we have
$$\widetilde{\psi}(\lambda,h,(\chi_{n})) \, = \, \lim_{n\rightarrow \infty}
- \frac{1}{n^{d-1}h(n)} \ln{ \mathbb P
  \left[\phi^{2^{\chi_{n}}}_{n^{d-1},h(n)} \geq \lambda n^{d-1} \right] }
\, = \, +\infty \, = \, \widetilde{\psi} (\lambda) \,. $$

$\bullet \,\, \lambda < \mu$ : We take $N,n \in \mathbb N$ with $n \leq N$ and let $N=nm+r$ be the Euclidean algorithm. We consider two functions $h$, $\widetilde{h}: \mathbb N \rightarrow \mathbb N$, with $\lim_{n\rightarrow \infty} h(n) =\lim_{n\rightarrow \infty}\widetilde{h}(n) =+\infty$, and let $\widetilde{h}(N)=h(n)\widetilde{m} + \widetilde{r}$ be the Euclidean algorithm. We take $k \in \mathbb N$ which will be chosen later. We want to compare $\phi^{k}_{N^{d-1},\widetilde{h}(N)}$ and $\phi^{k}_{n^{d-1},h(n)}$.

\begin{figure}[ht!]
\centering
\begin{picture}(0,0)%
\epsfig{file=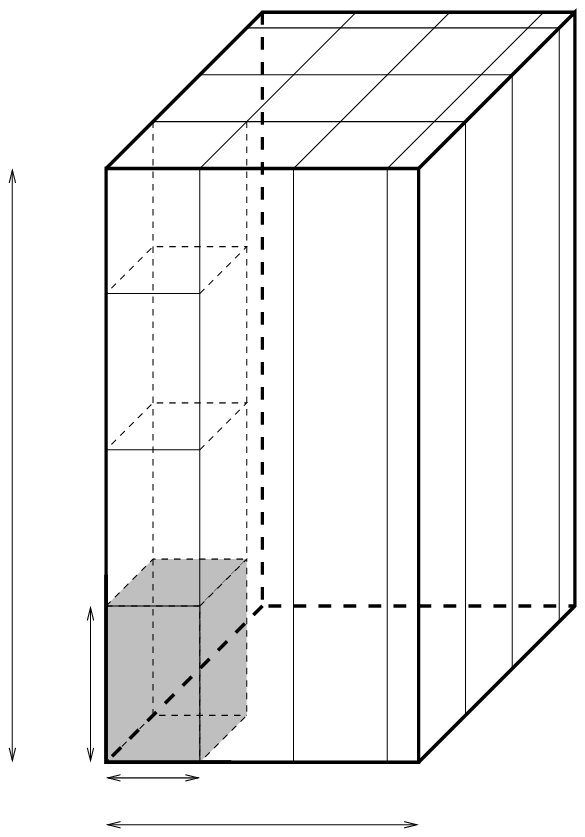}%
\end{picture}%
\setlength{\unitlength}{1973sp}%
\begingroup\makeatletter\ifx\SetFigFont\undefined%
\gdef\SetFigFont#1#2#3#4#5{%
  \reset@font\fontsize{#1}{#2pt}%
  \fontfamily{#3}\fontseries{#4}\fontshape{#5}%
  \selectfont}%
\fi\endgroup%
\begin{picture}(5583,8182)(2251,-8210)
\put(3751,-7636){\makebox(0,0)[b]{\smash{{\SetFigFont{6}{7.2}{\familydefault}{\mddefault}{\updefault}$n$}}}}
\put(4576,-8161){\makebox(0,0)[b]{\smash{{\SetFigFont{6}{7.2}{\familydefault}{\mddefault}{\updefault}$N$}}}}
\put(3076,-6586){\makebox(0,0)[rb]{\smash{{\SetFigFont{6}{7.2}{\familydefault}{\mddefault}{\updefault}$h(n)$}}}}
\put(2251,-4636){\makebox(0,0)[rb]{\smash{{\SetFigFont{6}{7.2}{\familydefault}{\mddefault}{\updefault}$\tilde{h}(N)$}}}}
\end{picture}%

\caption{Comparison between $\phi^{k}_{N^{d-1},\widetilde{h}(N)}$ and $\phi^{k}_{n^{d-1},h(n)}$}
\label{fig:cylindre}
\end{figure}

The idea is to divide $B((N,...,N),\widetilde{h}(N))$ into $m^{d-1}$ boxes which are disjoint translates of $B((n,...,n),\widetilde{h}(N))$, then to cut again $B((n,...,n),\widetilde{h}(N))$ into $\widetilde{m}$ disjoint translates of the elementary box $B((n,...,n),h(n))$ and to use here the lemma of junction (see figure \ref{fig:cylindre}).

We define two quantities that will allow us to deal with the edges
belonging to the part of  $\phi^{k}_{N^{d-1},\widetilde{h}(N)}$ that does not enter in any translate of $\phi^{k}_{n^{d-1},h(n)}$. On one hand, by the definition of $\mu$, $\lambda < \mu$ implies that $F([0,\lambda]) < 1$, so there exists a positive $\eta$ such that
$$ F([0,\lambda + \eta]) \, < \,1 , \,\,i.e., \,\, p(\eta)\, = \, \mathbb P [t(e) \geq \lambda + \eta] \, > \,0 \, .$$
It follows that there exists $k_{0}$ such that
$$\forall k\geq k_{0} \qquad \mathbb P \left[t^{k}(e) \geq \lambda + \frac{\eta}{2}\right] \, \geq \, p(\eta) \, > \,0 \, .  $$
On the other hand, if we define
$$ \gamma_{k} \, = \, \max \{ \, p\in \mathbb N \, | \, \mathbb P [t(e) \geq {p}{k}] >0 \, \} \wedge \left( \lfloor \lambda k n^{d-1} \rfloor +1 \right)  \, ,$$then we have
$$ p_{k} \, = \, \mathbb P \left[t(e) \geq \gamma_{k}k\right] \, > \,0 \, .$$
Let $k \geq k_{0}$ in $\mathbb N$. For $i_{1},...,i_{d-1}$ in $\{0,...,m \}$, we define
$$B_{i_{1},...,i_{d-1}} \, = \, \prod_{j=1}^{d-1} ]i_{j}n,(i_{j}+1)n] \times ]0,\widetilde{h}(N)] $$
and
$$B_{m^{d-1}+1} \, = \, B((N,...,N),\widetilde{h}(N)) \smallsetminus
\bigcup_{i_{1},...,i_{d-1} = 0}^{m-1} B_{i_{1},...,i_{d-1}}  \, .$$

\begin{rem}
\label{ajoutflux}
It is easy (and very useful) to see that if $C_{i} \times ]0,h]$, $i=1,2$
are two cylinders with disjoint bases $C_{1}$, $C_{2}
\subset \mathbb{R}^{d-1}$ having a common side and with maximal flows $\phi_{i}$, $i=1,2$, the
maximal flow through $(C_{1} \cup C_{2}) \times ]0,h]$ is at least
    $\phi_{1} + \phi_{2}$.
\end{rem}

We deduce from this remark that if for every $i_{1},...,i_{d-1}$ in $\{0,...,m-1\}$ we have $\phi^{k}_{B_{i_{1},...,i_{d-1}}} \geq \lambda n^{d-1}$ and if all the vertical edges $e$ in $B_{m^{d-1}+1}$ satisfy $t(e) \geq (\lambda +\eta)$, we have
$$ \phi^{k}_{N^{d-1},\widetilde{h}(N)} \, \geq \, \lambda N^{d-1} \,.$$
By independence we obtain
\begin{align*}
\mathbb P \left[ \phi^{k}_{N^{d-1},\widetilde{h}(N)} \geq \lambda N^{d-1} \right] & \, \geq \, \prod_{i_{1},...,i_{d-1} =0}^{m-1} \mathbb P \left[ \phi^{k}_{B_{i_{1},...,i_{d-1}}} \geq \lambda n^{d-1} \right] \times p(\eta) ^{(d-1)N^{d-2}r\widetilde{h}(N)} \\
& \, \geq \, \mathbb P \left[\phi^{k}_{n^{d-1},\widetilde{h}(N)} \geq \lambda n^{d-1}\right] ^{m^{d-1}} \times p(\eta) ^{(d-1)N^{d-2}r\widetilde{h}(N)} \, .\\
\end{align*}
We study next $\phi^{k}_{n^{d-1},\widetilde{h}(N)}$. We define for $j$ in $\{0,...,(\widetilde{m}-1)\}$
$$ B'_{j} \, = \, ]0,n]^{d-1} \times ]jh(n), (j+1)h(n)] $$
and
$$ B'_{\widetilde{m}} \, = \, B((n,...,n),\widetilde{h}(N)) \smallsetminus \bigcup_{j=0}^{\widetilde{m}-1} B'_{j} \, .$$
The probability of the boundary conditions $\Pi \in \{0,1/k,2/k,...,(
\lfloor \lambda n^{d-1} \rfloor +1 )/k\} ^{2n^{d-1}}$ in $B$ is the
probability that there exists a discrete stream $(g,o)$ in $B$ satisfying
$\Pi^{\lambda,n}(g)=\Pi$. Remember that every discrete stream $(g,o)$ must
satisfy the balance equation, so once we know that such a discrete stream
exists, $flow(g,o)$ is just given by the projection of $g$ on the vertical
edges that intersect the hyper-plane $\{(x_{1},...,x_{d}) \in
\mathbb{R}^{d} \, | \, x_{d} = h(n)-1/2  \}$, so $\Pi ^{\lambda ,n} (g)$
contains enough information to know if $flow(g,o)$ is bigger than $\lambda
n^{d-1}$ or not. We denote by $\Pi^{k}_{\lambda,n} = (\Pi^{k}_{\lambda,n,1}, \Pi^{k}_{\lambda,n,2})$ one of the boundary conditions of highest probability in $B((n,...,n),h(n))$ which corresponds to a discrete stream $(g,o)$ such that $ flow(g,o) \geq \lambda n^{d-1} $ and we define $(\Pi^{k}_{\lambda,n})^{*} = (\Pi^{k}_{\lambda,n,2}, \Pi^{k}_{\lambda,n,1})$. The model is invariant under reflections in the coordinates hyperplanes or translates of these hyperplanes, so by symmetry we have  $\mathbb P [\Pi^{k}_{\lambda,n}] = \mathbb P [ (\Pi^{k}_{\lambda,n})^{*}]$. Using the lemma of junction (lemma \ref{junction}), we know that if
\begin{itemize}
\item we can define a discrete stream in $B_{0}'$ with boundary conditions $\Pi^{k}_{\lambda,n}$,
\item we can define a discrete stream in $B_{1}'$ with boundary conditions $(\Pi^{k}_{\lambda,n})^{*}$,
\item we can define a discrete stream in $B_{2}'$ with boundary conditions $\Pi^{k}_{\lambda,n}$,
\item $\cdots$,
\item and all the vertical edges $e$ in $B_{\widetilde{m}}'$ satisfy $t(e)\geq \gamma_{k}k$,
\end{itemize}
then $\phi_{n^{d-1},\widetilde{h}(N)}^{k} \geq \lambda n^{d-1}$.
\begin{rem}
It is not sufficient to impose here that all the vertical edges $e$ in
$B'_{\widetilde{m}}$ satisfy $t(e) \geq \lambda$, because the amount of
fluid that goes out of $B'_{\widetilde{m}-1}$ at its top through one fixed edge $f$ 
can exceed $\lambda$ - we have no information about $\smash{\Pi^{k}_{\lambda ,
  n}}$ - and we cannot accept to lose fluid at the exit of $f$, unless it exceeds $\lambda n^{d-1}$. This is the reason why we introduced $\gamma_{k}$.
\end{rem}

Now by independence we obtain
\begin{equation}
\label{phi-Pi}
\mathbb P \left[ \phi^{k}_{n^{d-1},\widetilde{h}(N)} \geq \lambda
  n^{d-1}\right] \, \geq \, \mathbb P
\left[\Pi^{k}_{\lambda,n}\right]^{\widetilde{m}} \times p_{k}
^{n^{d-1}\widetilde{r}} \, ,
\end{equation}
whence
\begin{equation}
\label{lienN-n}
\mathbb P \left[\phi^{k}_{N^{d-1},\widetilde{h}(N)} \geq \lambda N^{d-1}\right] \, \geq \, \mathbb P \left[\Pi^{k}_{\lambda,n}\right]^{m^{d-1}\widetilde{m}} \times p_{k}^{n^{d-1}\widetilde{r} m^{d-1}} p(\eta)^{(d-1)N^{d-2}r\widetilde{h}(N)} \, .
\end{equation}
Let $\mathbf {\Pi}$ be the set of all the boundary conditions corresponding to a discrete stream $(g,o)$ such that $flow(g,o) \geq \lambda n^{d-1}$. We have seen that a maximal flow $\phi^{k}$ is always realized by a discrete stream, so we have
\begin{align*}
\mathbb P \left[\phi^{k}_{n^{d-1},h(n)} \geq \lambda n^{d-1}\right] & \, \leq \, \mathbb P \left[\bigcup_{\Pi \in \mathbf{\Pi}} \Pi \right] \\
& \, \leq \, \sum_{\Pi \in \mathbf{\Pi}} \mathbb P [\Pi] \\
& \, \leq \, N^{k}_{\lambda,n} \times \mathbb P \left[\Pi^{k}_{\lambda,n}\right] \, ,\\
\end{align*}
where we remember that $N^{k}_{\lambda,n}$ is the number of possible boundary conditions for discrete streams.

To obtain later a result independent of $k$, we need to consider two sequences $(k_{n}, n\in \mathbb N)$ and $(\widetilde{k}_{n}, n\in \mathbb N)$ such that $\lim_{n\rightarrow \infty} k_{n} = \lim_{n\rightarrow \infty} \widetilde{k}_{n} = +\infty$. We want to get rid of $N^{k}_{\lambda,n}$. We remember that
$$ N^{k_{n}}_{\lambda,n} \, \leq \, \left(k_{n}\left( \lfloor \lambda n^{d-1} \rfloor +1 \right) +1 \right)^{2n^{d-1}} \,.$$ 
Under the condition
\begin{equation}
\label{ln(kM)/h}
\lim_{n\rightarrow \infty} \frac{\ln{(k_{n}n)}}{h(n)} \, = \, 0
\end{equation}
we have
\begin{equation}
\label{lienphi-pi}
\limsup_{n\rightarrow \infty} \frac{1}{n^{d-1}h(n)} \ln{ \mathbb P \left[\phi^{k_{n}}_{n^{d-1},h(n)} \geq \lambda n^{d-1}\right] } \, \leq \, \limsup_{n\rightarrow \infty} \frac{1}{n^{d-1}h(n)} \ln{ \mathbb P \left[\Pi^{k_{n}}_{\lambda,n}\right] } \, .
\end{equation}

Consider (\ref{lienN-n}) again. This equation is satisfied for every $k
\geq k_{0}$, so it is true for $k_{n}$ with a fixed $n$ not too small. We
need to compare $t^{k_{n}}$ with $t^{\widetilde{k}_{N}}$, but the relation
is simple only if $\widetilde{k}_{N}$ is divisible by $k_{n}$. That is the
reason why from now on we will consider only sequences $(k_{n}, n\in \mathbb{N})$ and $(\widetilde{k}_{n}, n\in \mathbb N)$ such that
$$ \forall n\in \mathbb{N} \qquad k_{n} \,=\, 2^{\chi_{n}} \qquad and \qquad  \widetilde{k}_{n} \,=\, 2^{\widetilde{\chi}_{n}} $$
where $(\chi_{n},n\in \mathbb{N})$ and $(\widetilde{\chi}_{n},n\in
\mathbb{N})$ are non-decreasing sequences of integers. Of course the
condition $\lim_{n\rightarrow +\infty} k_{n} = \lim _{n\rightarrow +\infty}
\widetilde{k}_{n} = +\infty$ implies $\lim
_{n\rightarrow +\infty} \chi_{n} = \lim _{n\rightarrow +\infty}
\widetilde{\chi}_{n} = +\infty$. In that case for large $N$ we have
$\widetilde{\chi}_{N} \geq \chi_{n}$, so $\widetilde{k}_{N}$ is divisible
by $k_{n}$ and then $t^{\widetilde{k}_{N}} \geq t^{k_{n}}$, whence
\begin{equation}
\label{lienk,M}
\phi_{N^{d-1},\widetilde{h}(N)}^{\widetilde{k}_{N}} \geq \phi_{N^{d-1},\widetilde{h}(N)}^{k_{n}} \,.
\end{equation}
We use (\ref{lienN-n}) with $k=k_{n}=2^{\chi_{n}}$ and (\ref{lienk,M}) to obtain for $n$ and $N$ large enough
\begin{align*}
\frac{1}{N^{d-1}\widetilde{h}(N)} & \ln { \mathbb P\left[ \phi^{\widetilde{k}_{N}}_{N^{d-1},\widetilde{h}(N)} \geq \lambda N^{d-1} \right] } \\
& \, \geq \, \frac{1}{N^{d-1}\widetilde{h}(N)}  \ln { \mathbb P\left[ \phi^{k_{n}}_{N^{d-1},\widetilde{h}(N)} \geq \lambda N^{d-1} \right] } \\
& \, \geq \, \frac{m^{d-1}\widetilde{m}}{N^{d-1}\widetilde{h}(N)} \ln{\mathbb P \left[\Pi^{k_{n}}_{\lambda,n}\right]} + \frac{n^{d-1}m^{d-1}\widetilde{r}}{N^{d-1}\widetilde{h}(N)} \ln {p_{k_{n}}} + \frac{(d-1)N^{d-2}r\widetilde{h}(N)}{N^{d-1}\widetilde{h}(N)} \ln{p(\eta)} \, .\\
\end{align*}
We send first $N$ to $+\infty$ and then $n$ to $+\infty$; this gives us with the help of (\ref{lienphi-pi})
\begin{align*}
\liminf_{N\rightarrow \infty} \frac{1}{N^{d-1}\widetilde{h}(N)} & \ln { \mathbb P\left[ \phi^{\widetilde{k}_{N}}_{N^{d-1},\widetilde{h}(N)} \geq \lambda N^{d-1} \right] } \\
& \, \geq \, \limsup_{n\rightarrow \infty} \frac{1}{n^{d-1}h(n)} \ln{\mathbb P\left[\Pi^{k_{n}}_{\lambda,n}\right]} \\
& \, \geq \, \limsup_{n\rightarrow \infty} \frac{1}{n^{d-1}h(n)}  \ln { \mathbb P\left[ \phi^{k_{n}}_{n^{d-1},h(n)} \geq \lambda n^{d-1} \right] } \, .\\
\end{align*}
By considering the case $h=\widetilde{h}$ and
$k_{n}=\widetilde{k}_{n}=2^{\chi_{n}}$, under the condition
(\ref{ln(kM)/h}) on $h$ and $(k_{n})$ - i.e., the condition (\ref{h-chi}) on
$h$ and $(\chi_{n})$ -, we obtain the existence of the limit
$$ \widetilde{\psi}(\lambda,h,(\chi_{n})) \, = \, \lim_{n\rightarrow \infty} - \frac{1}{n^{d-1}h(n)}  \ln { \mathbb P\left[ \phi^{2^{\chi_{n}}}_{n^{d-1},h(n)} \geq \lambda n^{d-1} \right] } \, .  $$
For general $h$, $\widetilde{h}$, $\chi$ and $\widetilde{\chi}$ we obtain
that $\widetilde{\psi}(\lambda)$ is independent of the pair
$(h,(\chi_{n}) )$ satisfying (\ref{h-chi}), so theorem \ref{thmdiscret} is
proved.

\begin{rem}
\label{independancedepsi}
Thanks to this independence, we can prove some properties of
$\widetilde{\psi}$ by studying the behavior of the limit involved in
theorem \ref{thmdiscret} for specific choices of pairs $(h,(\chi_{n}))$.
\end{rem}

Moreover, still in the case $\lambda < \mu$, we have immediately that for $n$ sufficiently large
\begin{align*}
\mathbb P \left[\phi^{k_{n}}_{n^{d-1},h(n)} \geq \lambda n^{d-1}\right]
& \, \geq \, \mathbb P \left[\,all \,\, the \,\, vertical \,\, edges \,\, e \,\, in \,\, B((n,...,n),h(n)) \,\, satisfy \,\, t(e) \geq (\lambda + \eta)\, \right] \\
& \, \geq \, p(\eta)^{n^{d-1}h(n)} \, ,\\
\end{align*}
thus
$$ \widetilde{\psi}(\lambda) \, \leq \, -\ln{p(\eta)} \, < \, +\infty \, .$$

If the capacity $t$ of an edge is bounded by a constant $M$, we can simply define
$$ \forall x \in \mathbb{Z}^{d-1} \cap\, ]0,n]^{d-1} \qquad \pi_{h}^{\lambda,n}(g,x) \,=\, g\left(\langle (x,h),(x,h+1) \rangle\right) $$
without truncating $g$ because $g$ is already bounded by $M$. Then the number of possible boundary conditions $N^{k}_{n}$ satisfies
$$ N^{k}_{n} \,\leq\, \left( k(M+1) \right) ^{2n^{d-1}} $$
so we can replace the hypothesis (\ref{ln(kM)/h}) by
\begin{equation}
\label{borne}
\lim_{n\rightarrow \infty} \frac{\ln k_{n}}{h(n)} \,=\, 0 \,.
\end{equation}

\begin{rem}
We don't study the case $\lambda = \mu$ for the moment, it is more adapted to study it with the continuity of $\widetilde{\psi}$.
\end{rem}


\subsection{Convexity of $\widetilde{\psi}$}
\label{convexite}

Let $\lambda_{1} \leq \lambda_{2} < \mu$, and $\alpha \in ]0,1[$. We want to show that
\begin{equation}
\label{convex}
\widetilde{\psi}\left(\alpha \lambda_{1} + (1-\alpha) \lambda_{2} \right) \, \leq \, \alpha \widetilde{\psi}(\lambda_{1}) + (1-\alpha) \widetilde{\psi} (\lambda_{2}) \, .
\end{equation}
We know that $\widetilde{\psi}$ does not depend on the couple
$(h,(\chi_{n}))$ satisfying (\ref{h-chi}), so we can take $h(n)=n$ to
simplify the notations and we will take an adapted $(\chi_{n})$. First we fix $k$ in $\mathbb N$, we will make it vary later. We fix $n$, $m$ in $\mathbb N$, and take $N=nm$. We set $u=\lfloor \alpha m^{d-1} \rfloor$. We keep the same notations as in the previous section for $B_{i_{1},...,i_{d-1}}$, $i_{1},...,i_{d-1}$ in $\{0,...,m-1\}$. We use the lexicographic order to order $\{ (i_{1},...,i_{d-1}), \, i_{j} \in \{0,...,m-1\}, \, 1\leq j \leq (d-1) \}$ and use this to rename these cylinders $(B_{j}, 1\leq j \leq m^{d-1})$. On the event
$$ \left\{ \, \forall j \in \{1,...,u\}, \,\, \phi^{k}_{B_{j}} \geq \lambda_{1}\, \right\} \cap  \left\{ \, \forall j \in \{(u+1),...,m^{d-1}\}, \,\, \phi^{k}_{B_{j}} \geq \lambda_{2}\, \right\} $$
we have (see remark \ref{ajoutflux})
\begin{align*}
\phi^{k}_{N^{d-1},N} & \, \geq \, \left( u\lambda_{1}n^{d-1} + (m^{d-1}-u)\lambda_{2}n^{d-1} \right) \\
& \, \geq \, N^{d-1} \left( \frac{u}{m^{d-1}} \lambda_{1} + \left(1-\frac{u}{m^{d-1}}\right) \lambda_{2} \right) \\
& \, \geq \, N^{d-1} \left( \alpha \lambda_{1} + (1-\alpha)\lambda_{2} \right) \\
\end{align*}
because $\lambda_{1} < \lambda_{2}$, so
$$ \mathbb P \left[\phi^{k}_{N^{d-1},N} \geq N^{d-1} \left(\alpha \lambda_{1}+(1-\alpha)\lambda_{2}\right) \right] \, \geq \, \mathbb P \left[\phi^{k}_{n^{d-1},N} \geq \lambda_{1}n^{d-1}\right] ^{u} \times \mathbb P \left[\phi^{k}_{n^{d-1},N} \geq \lambda_{2}n^{d-1}\right] ^{m^{d-1}-u} \, .$$
As in the previous section (see (\ref{phi-Pi})), we have
$$ \mathbb P \left[\phi^{k}_{n^{d-1},N} \geq \lambda_{i} n^{d-1}\right]  \,
\geq \, \mathbb P \left[\Pi^{k}_{\lambda_{i},n}\right]^{m}  \qquad i=1,2 \,.$$
so
\begin{align*}
\frac{1}{N^{d}} \ln \mathbb{P} \left[ \right. & \left. \phi^{k}_{N^{d-1},N}  \geq N^{d-1} \left(\alpha \lambda_{1}+(1-\alpha)\lambda_{2}\right)\right] \\ & \geq \, \frac{m^{d-1}}{N^{d}} \left( \frac{u}{m^{d-1}} \ln { \mathbb P \left[\Pi^{k}_{\lambda_{1},n}\right]} + \left( 1-\frac{u}{m^{d-1}}\right) \ln { \mathbb P \left[\Pi^{k}_{\lambda_{2},n}\right]}  \right) \, .\\
\end{align*}
We make now $k$ vary, $k_{n}=2^{\chi_{n}}$ with $(n,(\chi_{n}))$ satisfying
the condition (\ref{h-chi}) (for example $\chi_{n}=\lfloor n^{1/2}
\rfloor$), and we use the property $\phi_{N^{d-1},N}^{2^{\chi_{N}}} \geq
\phi_{N^{d-1},N}^{2^{\chi_{n}}}$ for large $N$; we send first $N$ to $+\infty$ and then $n$ to $+\infty$. We proved in the previous section that
$$ \limsup_{n\rightarrow \infty} \frac{1}{n^{d-1}h(n)} \ln \mathbb P \left[ \Pi_{\lambda,n}^{2^{\chi_{n}}} \right] \,=\, - \widetilde{\psi}(\lambda) \,, $$
so we obtain (\ref{convex}).


\subsection{Continuity of $\widetilde{\psi}$}
\label{continuite}

We want to show that $\widetilde{\psi}$ is continuous on $[0,\mu]$ when
$\mu$ is finite or on $[0,+\infty[$ when $\mu$ is infinite (remember that $\widetilde{\psi}$ is infinite on $]\mu,
  +\infty[$). The function $\widetilde{\psi}$ is convex and finite
    on $[0,\mu[$, so $\widetilde{\psi}$ is continuous on
        $]0,\mu[$. We assume then that
        $0 < \mu$: $\widetilde{\psi}(0)=0$ and $\widetilde{\psi}$
        is non-negative on $\mathbb R^{+}$, so $\widetilde{\psi}$ is right
        continuous at $0$. We assume then that
        $0 < \mu< +\infty$. The only point which remains to study is the
        left continuity of $\widetilde{\psi}$ at $\mu$. Remember that we
        did not define $\widetilde{\psi}$ at $\mu$, we will do it now. We set
$$ q_{\mu} \, = \, \mathbb P [t(e)=\mu] \,. $$
Notice that $q_{\mu}$ can be null. We remark that
\begin{align*}
\mathbb P \left[\phi_{n^{d-1},h(n)} \geq \mu n^{d-1}\right] & \, = \, \mathbb P \left[\, all \,\, the \,\, vertical \,\, edges \,\, e \,\, in \,\, B((n,...,n),h(n)) \,\, satisfy \,\, t(e)=\mu \, \right] \\
& \, = \, q_{\mu}^{n^{d-1}h(n)} \, ,\\
\end{align*}
so
$$ \lim_{n\rightarrow \infty} -\frac {1}{n^{d-1}h(n)} \ln { \mathbb P
  \left[\phi_{n^{d-1},h(n)} \geq \mu n^{d-1} \right] } \, = \, -\ln{q_{\mu}}$$
is finite as soon as $q_{\mu} >0$. Unfortunately, the existence of an atom for the law of $t(e)$ at $\mu$ does not imply the existence of an atom for the law of $t^{k}(e)$ at $\mu$, so we can have $q_{\mu} >0$ and
$$  \lim_{n\rightarrow \infty} -\frac {1}{n^{d-1}h(n)} \ln { \mathbb P
  \left[\phi^{k_{n}}_{n^{d-1},h(n)} \geq \mu n^{d-1} \right] } \, = \, + \infty \, .$$
This is the reason why we did not study $\widetilde{\psi}(\mu)$
  previously. We define (for every pair $(h,(\chi_{n}))$ as in theorem \ref{thmdiscret})
$$ \widetilde{\psi}(\mu) \, = \, - \ln{q_{\mu}} \, ,$$
which can eventually be infinite.

Now we want to check that $\widetilde{\psi}$ is left continuous at
$\mu$ (if $q_{\mu}=0$ we will show that $\lim \widetilde{\psi} (\lambda) =
+\infty$ when $\lambda\leq \mu$ and $\lambda \rightarrow \mu$). The idea of the proof is simple: if the flow in a cylinder is big,
it must be big in each horizontal section of this cylinder. We fix
$\varepsilon >0$, and we take $h(n),k_{n} \rightarrow_{n \rightarrow +\infty} +\infty$, $k_{n}=2^{\chi_{n}}$, satisfying the condition (\ref{ln(kM)/h}). We define for $i$ in $\{0,...,(h(n)-1) \}$
$$ C_{i} \, = \, ]0,n]^{d-1} \times ]i,i+1] $$
and we denote by $t_{1},...,t_{n^{d-1}}$ the capacities of the $n^{d-1}$ vertical edges in $C_{0}$. We have
\begin{align*}
\mathbb P \left[\phi^{k_{n}}_{n^{d-1},h(n)} \geq (\mu -\varepsilon) n^{d-1} \right] & \, \leq \, \mathbb P \left[ \bigcap_{i=0}^{h(n)-1} \{\, \phi^{k_{n}}_{C_{i}} \geq (\mu - \varepsilon) n^{d-1}\, \} \right] \\
& \, \leq \, \mathbb P \left[\phi^{k_{n}}_{n^{d-1},1} \geq (\mu - \varepsilon) n^{d-1} \right]^{h(n)} \,,\\
\end{align*}
and we know that
$$\phi_{n^{d-1},1}^{k_{n}} \, =\, \sum _{j=1}^{n^{d-1}} t_{j}^{k_{n}} \,
\leq \, \sum_{j=1}^{n^{d-1}} t_{j} \,,$$
so we have
$$\mathbb P \left[\phi^{k_{n}}_{n^{d-1},h(n)} \geq (\mu -\varepsilon)
  n^{d-1} \right]  \, \leq \, \mathbb{P} \left[ \sum_{j=1}^{n^{d-1}}
  (t_{j}-\mu) \geq -\varepsilon n^{d-1}  \right]^{h(n)} \,. $$
For every positive $\rho$ we obtain
$$\mathbb P \left[\phi^{k_{n}}_{n^{d-1},h(n)} \geq (\mu -\varepsilon) n^{d-1} \right]  \, \leq \, e^{\rho \varepsilon n^{d-1} h(n)} \mathbb E [e^{\rho (t-\mu)}]^{n^{d-1}h(n)} \, .$$
This expectation is well defined, because $(t-\mu)\leq 0$. Let $\eta >0$. Since
$$ \lim_{\rho \rightarrow +\infty} \mathbb E[e^{\rho(t-\mu)}] \, = \, q_{\mu} \, ,$$
then there exists $\rho_{0}$ such that
$$ \forall \rho \geq \rho_{0} \qquad  \mathbb E[e^{\rho(t-\mu)}]  \, \leq \, (q_{\mu}+\eta ) \, .$$
It follows that
$$ \frac{1}{n^{d-1}h(n)} \ln {\mathbb P \left[\phi^{k_{n}}_{n^{d-1},h(n)} \geq (\mu -\varepsilon) n^{d-1}\right]} \, \leq \, \rho_{0}\varepsilon + \ln{(q_{\mu}+\eta)} \, ,$$
so
$$ \widetilde{\psi}(\mu-\varepsilon) \, \geq \, -\rho_{0}\varepsilon - \ln{(q_{\mu}+\eta)} \, ,$$
whence
$$ \lim_{\varepsilon \rightarrow 0} \widetilde{\psi}(\mu-\varepsilon) \, \geq \, - \ln{(q_{\mu}+\eta)} \, .$$
This is true for every positive $\eta$, so
$$ \lim_{\varepsilon \rightarrow 0} \widetilde{\psi}(\mu-\varepsilon) \, \geq \, \lim_{\eta \rightarrow 0} - \ln{(q_{\mu}+\eta)} \, = \, -\ln{q_{\mu}} \, = \, \widetilde{\psi} (\mu) \, .$$
If $q_{\mu} =0$, we have the desired equality. Otherwise, we remark that
for every positive $\varepsilon$ we have
\begin{align*}
\mathbb P \left[ \phi^{k_{n}}_{n^{d-1},h(n)} \geq \right. & \left. (\mu-\varepsilon) n^{d-1} \right] \\
& \,\geq\, \mathbb P \left[ all\,\, the\,\,vertical\,\,edges\,\,e\,\,in\,\,B((n,...,n),h(n))\,\,satisfy\,\,t^{k_{n}}(e)\geq(\mu-\varepsilon) \right] \\
& \,\geq\, \mathbb P \left[ t^{k_{n}}(e) \geq \mu-\varepsilon \right] ^{n^{d-1}h(n)} \,. \\
\end{align*}
Now for $k_{n}$ sufficiently large we have
$$ \mathbb P \left[ t^{k_{n}}(e) \geq \mu-\varepsilon \right] \,\geq\, \mathbb P \left[ t(e)\geq \mu-\frac{\varepsilon}{2} \right] \,\geq\, q_{\mu} \,, $$
thus
$$ \forall \varepsilon>0 \qquad \widetilde{\psi}(\mu-\varepsilon) \,\leq\, -\ln q_{\mu} \,=\, \widetilde{\psi}(\mu) \,.$$
This ends the proof of the continuity of $\widetilde{\psi}$ on $[0,\mu]$
(or $[0,+\infty [$ if $\mu$ is infinite). We deduce immediately from this
    continuity that $\widetilde{\psi}$ is good.


\subsection{Existence of the limit for $\phi_{n^{d-1},h(n)}$}
\label{limitegen}

We come back to the existence of the limit involving $\phi$ in theorem \ref{limite}. We consider three cases.

$\bullet \,\, \lambda > \mu$ : Then
$$\forall n\in \mathbb N \qquad \mathbb P \left[\phi_{n^{d-1},h(n)} \geq \lambda n^{d-1}\right] \, = \, 0 \, ,$$
so the limit involved in theorem \ref{limite} exists and satisfies
$$ \psi(\lambda) \,=\, \lim_{n\rightarrow \infty} -\frac {1}{n^{d-1}h(n)} \ln{\mathbb P\left[\phi_{n^{d-1},h(n)} \geq \lambda n^{d-1}\right]} \, = \, + \infty \, = \, \widetilde{\psi}(\lambda) \,. $$

$\bullet \,\, \lambda = \mu$ : As we saw by studying the continuity of $\widetilde{\psi}$, we have
$$ \psi(\mu) \,=\, \lim_{n\rightarrow \infty} -\frac {1}{n^{d-1}h(n)} \ln{\mathbb P\left[\phi_{n^{d-1},h(n)} \geq \mu n^{d-1}\right]} \, = \, -\ln{q_{\mu}} \, = \, \widetilde{\psi}(\mu) $$
by definition of $\widetilde{\psi}(\mu)$.

$\bullet \,\, \lambda <\mu $ : We will compare $\phi^{k_{n}}_{n^{d-1},h(n)}$ with $\phi_{n^{d-1},h(n)}$. We fix $k \in \mathbb{N}$. We know that $t(e) \geq t^{k}(e)$ so $\phi_{n^{d-1},h(n)} \geq \phi^{k}_{n^{d-1},h(n)} $. For a set of edges $E$, we denote by $V^{k}(E)$ the quantity $\sum_{e\in E} t^{k}(e)$. Thanks to the max-flow min-cut theorem we obtain
$$ \phi^{k}_{n^{d-1},h(n)} \, = \, \min \{\, V^{k}(E) \, | \, E \,\, is \,\, an \,\, (F_{0},F_{h(n)})-cut\,\}\,. $$
Let $E_{0}$ be an $(F_{0},F_{h(n)})$-cut realizing this minimum (it may
depend on $k$). Then
\begin{align*}
\phi^{k}_{n^{d-1},h(n)} & \, = \, V^{k}(E_{0}) \\
& \, = \, \sum_{e\in E_{0}} t^{k}(e) \\
& \, \geq \, \sum_{e\in E_{0}} t(e) - \frac{|E_{0}|}{k}  \\
& \, \geq \, \min \{\, V(E) \, | \, E \,\, is \,\, an \,\, (F_{0},F_{h(n)})-cut\, \} - n^{d-1} \frac{h(n)}{k} \\
& \, \geq \, \phi_{n^{d-1},h(n)} - n^{d-1}\frac{h(n)}{k} \, .\\
\end{align*}
We fix $\lambda \geq 0$, and we make now $k$ vary. We take $k_{n}=
2^{\chi_{n}}$ (with $(\chi_{n})$ a non-decreasing sequence of integers such
that $\lim_{n\rightarrow +\infty} \chi_{n} = +\infty$) satisfying with $h$ the condition (\ref{ln(kM)/h}). If the sequence $(k_{n})$ satisfies also the condition
\begin{equation}
\label{h/k}
\lim_{n\rightarrow \infty} \frac{h(n)}{k_{n}} \,=\,0
\end{equation}
then we have for every $\lambda' <\lambda$ the existence of $n_{0} \in \mathbb{N}$ such that
$$\forall n\geq n_{0} \qquad \lambda - \frac{h(n)}{k_{n}} \,\geq\, \lambda' \,. $$
We deduce that under the condition (\ref{h/k}) we have for all $n\geq n_{0}$
$$ \mathbb P \left[ \frac{\phi^{k_{n}}_{n^{d-1},h(n)}}{n^{d-1}} \geq \lambda \right] \,\leq\,  \mathbb P \left[ \frac{\phi_{n^{d-1},h(n)}}{n^{d-1}} \geq \lambda \right] \,\leq\,  \mathbb P \left[ \frac{\phi^{k_{n}}_{n^{d-1},h(n)}}{n^{d-1}} \geq \lambda' \right] \,. $$
We conclude thanks to the hypothesis (\ref{ln(kM)/h}) that
$$ \widetilde{\psi}(\lambda) \, \geq \, \limsup_{n\rightarrow \infty}\, (\square) \, \geq \, \liminf_{n\rightarrow \infty}\, (\square) \, \geq \, \widetilde{\psi}(\lambda') $$
where
$$ (\square)\, = \,-\frac{1}{n^{d-1}h(n)} \ln {\mathbb P \left[\phi_{n^{d-1},h(n)} \geq \lambda n^{d-1}\right]} \, .$$
Sending $\lambda'$ to $\lambda$, thanks to the continuity of $\widetilde{\psi}$ in $[0,\mu[$, we obtain the existence of the limit
$$ \psi(\lambda, h) \,=\, \lim_{n\rightarrow \infty} -\frac{1}{n^{d-1}h(n)}
    \ln{\mathbb P\left[\phi_{n^{d-1},h(n)} \geq \lambda n^{d-1}\right]} \,
    = \, \widetilde{\psi}(\lambda) \, .$$
Moreover we know that this limit is independent of $h$ satisfying
    $\lim_{n\rightarrow +\infty} h(n) = +\infty$ and such that there exists a
    non-decreasing sequence of integers $(\chi_{n})$, $\lim_{n\rightarrow
    +\infty} \chi_{n} = +\infty$ for which the pair $(h, (2^{\chi_{n}}))$
    satisfies (\ref{ln(kM)/h}) and (\ref{h/k}): we denote it by
    $\psi(\lambda)$. It is finally obvious that the existence of such a sequence $(\chi_{n})$ is equivalent to the condition
$$ \lim_{n\rightarrow \infty} \frac{h(n)}{\ln n} \,=\, +\infty $$
(let $\chi_{n}= \lfloor 2 \ln h(n) /\ln 2 \rfloor $ for example). This ends
    the proof of the existence of the limit $\psi$ in theorem \ref{limite},
    and we have $\psi=\widetilde{\psi}$ so the properties proved for
    $\widetilde{\psi}$ still hold for $\psi$.

If the capacity $t$ of an edge is bounded, we can replace the condition (\ref{ln(kM)/h}) by (\ref{borne}); in that case, as soon as
$$ \lim_{n\rightarrow \infty} h(n) \,=\,+\infty $$
we can find a sequence $(k_{n})=(2^{\chi_{n}})$ satisfying (\ref{borne}) and (\ref{h/k}), so the limit exists.


\subsection{The function $\psi $ vanishes on $[0,\nu (F)]$}

This could be proved easily thanks to theorem \ref{lgn} in dimension three
and with the hypothesis on $F$ required in theorem \ref{lgn}, but we prefer to prove it directly in the general case without
theorem \ref{lgn}.

\begin{figure}[ht!]
\centering

\begin{picture}(0,0)%
\epsfig{file=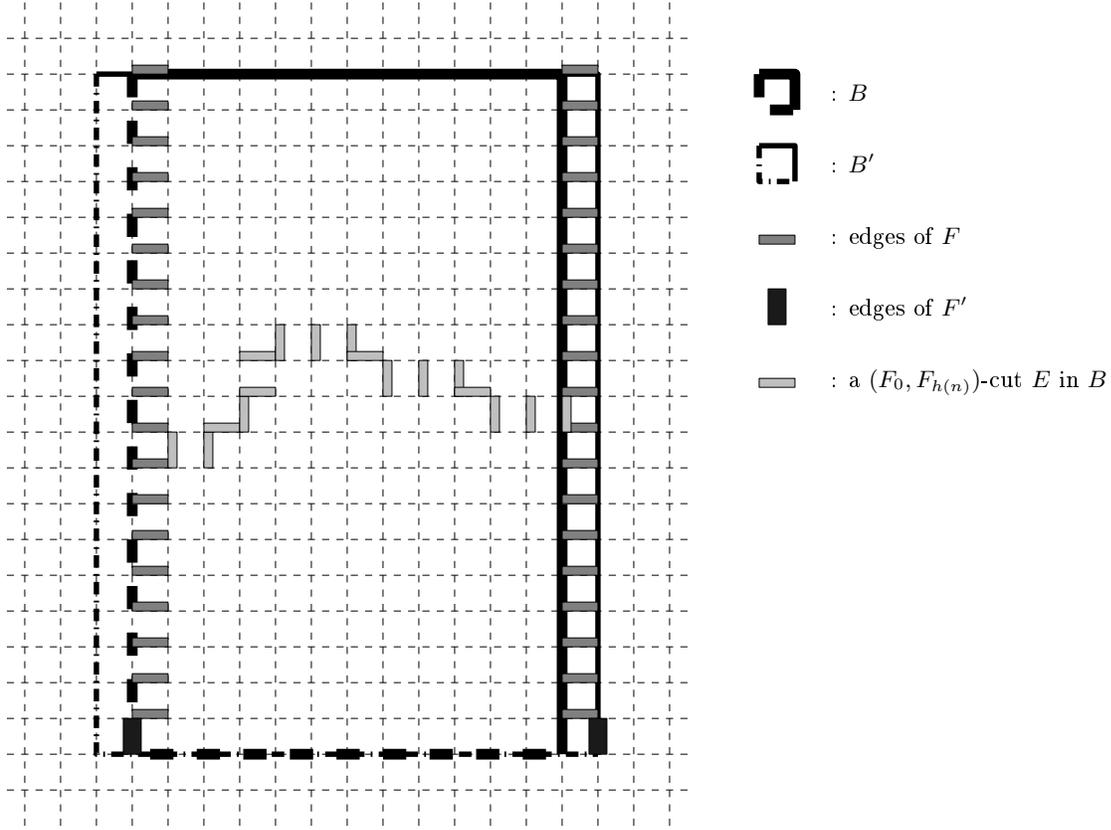}%
\end{picture}%
\setlength{\unitlength}{2960sp}%
\begingroup\makeatletter\ifx\SetFigFont\undefined%
\gdef\SetFigFont#1#2#3#4#5{%
  \reset@font\fontsize{#1}{#2pt}%
  \fontfamily{#3}\fontseries{#4}\fontshape{#5}%
  \selectfont}%
\fi\endgroup%
\begin{picture}(9927,6924)(2539,-8473)
\put(9451,-2386){\makebox(0,0)[lb]{\smash{{\SetFigFont{9}{10.8}{\rmdefault}{\mddefault}{\updefault}{\color[rgb]{0,0,0}: $B$}%
}}}}
\put(9451,-2986){\makebox(0,0)[lb]{\smash{{\SetFigFont{9}{10.8}{\rmdefault}{\mddefault}{\updefault}{\color[rgb]{0,0,0}: $B'$}%
}}}}
\put(9451,-3586){\makebox(0,0)[lb]{\smash{{\SetFigFont{9}{10.8}{\rmdefault}{\mddefault}{\updefault}{\color[rgb]{0,0,0}: edges of $F$}%
}}}}
\put(9451,-4186){\makebox(0,0)[lb]{\smash{{\SetFigFont{9}{10.8}{\rmdefault}{\mddefault}{\updefault}{\color[rgb]{0,0,0}: edges of $F'$}%
}}}}
\put(9451,-4786){\makebox(0,0)[lb]{\smash{{\SetFigFont{9}{10.8}{\rmdefault}{\mddefault}{\updefault}{\color[rgb]{0,0,0}: a $(F_{0},F_{h(n)})$-cut $E$ in $B$}%
}}}}
\end{picture}%

\caption{Comparison between $\phi$ and $\tau$ in dimension two}
\label{phitau}
\end{figure}

We suppose now that $\mathbb{E} [t]$ is finite. We suppose that $\nu >0$ (otherwise there is nothing to prove),
and we take $\lambda = \nu - \varepsilon$, with a positive
$\varepsilon$. Remark \ref{independancedepsi} holds for $\psi$ too: we know
that $\psi$ is independent of $h$ satisfying $\lim _{n\rightarrow +\infty}
h(n) / \ln n = +\infty$ so we can make a specific choice of function $h$ and
study the corresponding limit to show a general result on $\psi$. We take $h\rightarrow \infty$ such that
\begin{equation}
\label{hpetit}
\lim_{n\rightarrow \infty} \frac{h(n)}{n}\, = \, 0   \qquad  and \qquad \lim_{n\rightarrow \infty} \frac{h(n)}{\ln n} \,=\, +\infty  \,.
\end{equation}
We remember that
$$ \tau_{n^{d-1}} \, = \, \tau (]0,n]^{d-1}) \, = \, \inf\, \{\, V(E) \, |\, E \,\, is \,\, a \,\, cut \,\, over \,\, ]0,n]^{d-1} \,\, and \,\, E \,\, satisfies \,\, (*)\, \} \, ,$$
where $(*)$ is defined at the end of the section \ref{*}. We define for $S$ a hyper-rectangle the variable
$$  \tau (S,k) \, = \, \inf\, \{\, V(E) \, |\, E \,\, is \,\, a \,\, cut \,\, over \,\, S \,\, , \,\, E \,\, satisfies \,\, (*) \,\, and \,\, E\subset S\times]-k,k]\, \} \, ,$$
and
$$ \tau_{n^{d-1},k}\, = \, \tau (]0,n]^{d-1},k) \, .$$
We define the set of edges $F$ as
$$ F \, =\, \{ \langle x,y \rangle \,|\, x \in B \,,\,\, y \notin B \,\,
and \,\, \langle x,y\rangle \in \mathbb{R}^{d-1}\times [1,h(n)] \} \, . $$
This is the set of the edges through which some fluid could escape from $B$
somewhere else than at its bottom or at its top. We denote by $|F|$ the cardinality of $F$, $|F|= 2(d-1)n^{d-2}h(n)$. We consider the larger cylinder
$$B'\,=\,  ]-1,n+1]^{d-1} \times ]0,h(n)] \,, $$
and we define
$$ \tau'_{(n+2)^{d-1},h(n)} \,=\, \tau \left( ]-1,n+1]^{d-1},h(n) \right) \,. $$
We finally define the set of edges
$$ F' \,=\, \{\, e\in B'\smallsetminus B \,|\,e\,\,is\,\,vertical\,, \,\, e\in \mathbb{R}^{d-1}\times [0,1] \} \, $$
of cardinality $|F'|=2(d-1)(n+1)^{d-2}$ (see figure \ref{phitau} in
dimension two). 
We remark that if $E$ is an $(F_{0},F_{h(n)})$-cut in $B((n,...,n),h(n))$,
the set of edges $E\cup F \cup F' $ contains a cut over $
]-1,n+1]^{d-1}$ satisfying the condition $(*)$ for $S=]-1,n+1]^{d-1}$, so
$$  \tau'_{(n+2)^{d-1},h(n)} - \phi_{n^{d-1},h(n)}  \, \leq \, \sum_{e\in F } t(e) + \sum_{e\in F'} t(e) \, .$$
We obtain for $M>\mathbb E[t]$
\begin{align*}
\mathbb P\left[\frac{\phi_{n^{d-1},h(n)}}{n^{d-1}}\geq\lambda \right] & \, \geq \, \mathbb P \left[\left\{\, \frac{ \phi_{n^{d-1},h(n)}}{n^{d-1}} \geq \lambda\, \right\} \cap \left\{\, \frac{ \tau'_{(n+2)^{d-1},h(n)}- \phi_{n^{d-1},h(n)} }{|F|+|F'|} \leq M\, \right\} \right] \\
& \, \geq \, \mathbb P \left[\left\{ \, \frac{ \tau'_{(n+2)^{d-1},h(n)}}{n^{d-1}} \geq \lambda +M \frac{|F|+|F'|}{n^{d-1}}\, \right\} \right. \\
& \qquad \qquad \left. \cap \left\{\, \frac{ \tau'_{(n+2)^{d-1},h(n)}-
    \phi_{n^{d-1},h(n)} }{|F|+|F'|} \leq M\, \right\} \right]\,. \\
\end{align*}
We remark that $\tau'_{(n+2)^{d-1},h(n)}$ is equal in law to
$\tau_{(n+2)^{d-1},h(n)}$, so
\begin{align*}
\mathbb P\left[\frac{\phi_{n^{d-1},h(n)}}{n^{d-1}}\geq\lambda \right] 
& \, \geq \, 1 - \Bigl( \mathbb P \left[\frac{ \tau_{(n+2)^{d-1},h(n)}}{n^{d-1}} < \lambda +M\frac{|F|+|F'|}{n^{d-1}}\right] \\
& \qquad \qquad + \mathbb P \left[ \frac{ \tau'_{(n+2)^{d-1},h(n)} -
    \phi_{n^{d-1},h(n)} }{|F|+|F'|} > M\right] \Bigr)  \\
& \, \geq \, 1- \Bigl( \mathbb P \left[\frac{ \tau_{(n+2)^{d-1}}}{n^{d-1}} < \nu-\frac{\varepsilon}{2}\right] + \mathbb P\left[ \frac{1}{|F|+|F'|} \sum_{e\in F \cup F'} t(e) \geq M \right] \Bigr) \\
\end{align*}
for $n$ sufficiently large, thanks to (\ref{hpetit}) and the fact that
$\tau_{(n+2)^{d-1},h(n)}\geq \tau_{(n+2)^{d-1}}$. We know that $M>\mathbb E[t]$ and $\lim_{n\rightarrow \infty} (\tau_{(n+2)^{d-1}} / n^{d-1}) = \nu$ almost surely, so
$$ \lim_{n\rightarrow \infty} \mathbb P\left[\frac{\phi_{n^{d-1},h(n)}}{n^{d-1}}\geq\lambda \right] \, = \, 1 \, ,$$
which leads to
$$\psi(\lambda)\, = \,0 \, .$$
To conclude that $\psi(\nu)=0$ we need only to check that $\psi$ is left
continuous at $\nu$, i.e., to be sure that $\nu \leq \mu$. Suppose that $\nu
> \mu$, then $\mathbb P[t\geq \nu] =0$, so $\mathbb E[t] < \nu$, and we can
find a positive $\varepsilon$ such that $\mathbb E [t] < \nu -
\varepsilon$. Now if we denote by $(\widetilde{t}_{i}, i=1,...,n^{d-1})$
the capacities of the vertical edges in $]0,n]^{d-1}\times ]0,1]$, we have
$$ \mathbb P \left[\frac{\tau_{n^{d-1}}}{n^{d-1}} \geq \nu-\varepsilon\right] \, \leq \, \mathbb P \left[\frac{\sum_{i=1}^{n^{d-1}} \widetilde{t}_{i}}{n^{d-1}} \geq \nu-\varepsilon\right] \xrightarrow[n\rightarrow \infty]{} 0 \, .$$
This is absurd because $(\tau_{n^{d-1}}/n^{d-1})$ converges toward $\nu$ almost surely. We conclude that $\nu \leq \mu$ and that $\psi(\nu) = 0$.


\subsection{The function $\psi $ is positive on $]\nu (F), +\infty[$}

We suppose that there exists a positive $\theta$ such that $\int_{[0,+\infty[} e^{\theta x} dF(x)$ is finite. The proof is based on the Cram\'er theorem in $\mathbb R$.

Let $\lambda = \nu + \varepsilon$, for a positive $\varepsilon$. We fix $k$, $N \in \mathbb N$, we will choose them later. We define
$$u \, = \, \lfloor \frac{h(N)}{2k} \rfloor \, .$$
Just as in the study of the continuity of $\psi$, by cutting  $B((N,...,N),h(N))$ into horizontal sections of height $2k$, we have
$$ \mathbb P \left[\phi_{N^{d-1},h(N)} \geq (\nu+\varepsilon)N^{d-1}\right]
\, \leq \, \mathbb P \left[ \phi_{N^{d-1},2k} \geq (\nu+\varepsilon)
  N^{d-1}\right]^{u} \,=\, \mathbb{P} \left[ \phi_{\mathcal{B}(k)} \geq
  (\nu+\varepsilon) N^{d-1} \right]^{u} \, ,$$
where $\mathcal{B}(k)=]0,N^{d-1}]\times ]-k,k]$ because $\phi_{N^{d-1},2k}$
    and $\phi_{\mathcal{B}}$ are equal in law.
Now $\mathbb E [\tau(S,k)]$ is subadditive in the sense that for disjoint
hyper-rectangles $S_{1}$ and $S_{2}$ having a common side, we have
$$ \tau(S_{1}\cup S_{2},k) \, \leq \, \tau (S_{1},k) + \tau (S_{2},k) \, .$$
Moreover $\mathbb E [\tau(S,k)]$ is non-negative and finite (because $\mathbb E [t] <\infty$), so by a classical subadditive argument we have the existence of
$$ \nu_{k} \, = \, \lim_{n\rightarrow \infty} \frac{\mathbb E [\tau _{n^{d-1},k}]}{n^{d-1}} $$
and we know that
$$ \nu_{k}\, = \, \inf_{n} \frac{\mathbb E [\tau _{n^{d-1},k}]}{n^{d-1}} \, .$$
The sequence $(\nu_{k},k\in \mathbb N)$ is non-increasing in $k$ and non-negative, so it converges; we denote by $\widetilde{\nu}$ its limit: $\widetilde{\nu}=\lim_{k\rightarrow \infty} \nu_{k} = \inf_{k} \nu_{k}$.  By the same subadditive argument, we have
$$ \lim_{n\rightarrow \infty} \frac{\mathbb E[\tau_{n^{d-1}}]}{n^{d-1}} \, = \, \nu \, = \, \inf_{n} \frac{\mathbb E[\tau_{n^{d-1}}]}{n^{d-1}} \, .$$
We obtain
$$ \widetilde{\nu} \, = \, \inf_{k} \inf_{n} \frac{\mathbb E [\tau _{n^{d-1},k}]}{n^{d-1}} \, = \, \inf_{n} \inf_{k} \frac{\mathbb E [\tau _{n^{d-1},k}]}{n^{d-1}} \, = \, \nu \, ,$$
thus we can choose $k_{0}$ such that $\nu_{k_{0}} \leq \nu + \varepsilon/4$. Then we choose $n_{0}$ such that
$$ \frac{\mathbb E [\tau_{n_{0}^{d-1},k_{0}}]}{n_{0}^{d-1}} \, < \, \nu_{k_{0}} + \frac{\varepsilon}{2} \, ,$$
and we fix $N=n_{0}m$, with $m \in \mathbb N$. We have
$$ \phi_{\mathcal{B}(k_{0})} \, \leq \, \tau_{N^{d-1},k_{0}} \, \leq \, \sum_{i_{1},...,i_{d-1} = 0}^{m-1} \tau \left(\prod_{j=1}^{d-1} ]i_{j}n_{0}, (i_{j}+1)n_{0}],k_{0}\right) \, .$$
The variables $(\tau (\prod_{j=1}^{d-1} ]i_{j}n_{0}, (i_{j}+1)n_{0}],k_{0}), 0\leq i_{1},...,i_{d-1} \leq m-1)$ are independent and identically distributed, with the same law as $\tau_{n_{0}^{d-1},k_{0}}$. Their common expectation is
$$ \mathbb E [\tau_{n_{0}^{d-1},k_{0}}] \, \leq \, \left( \nu_{k_{0}} +
    \frac{\varepsilon}{2} \right) n_{0}^{d-1} \, .$$
Moreover for some positive $\theta$ we know that $\mathbb{E} [e^{\theta
    t}]$ is finite so
$$ \mathbb E \left[ e^{\theta \tau_{n_{0}^{d-1},k_{0}}} \right] \, \leq \,
\mathbb E \left[ e^{\theta \sum_{i=1}^{n_{0}^{d-1}} \widetilde{t}_{i}}\right] \, \leq \, \mathbb E \left[e^{\theta t}\right] ^{n_{0}^{d-1}} \, < \, \infty \, ,$$
where  $(\widetilde{t}_{i}, 1\leq i \leq n_{0}^{d-1})$ are still the
capacities of the vertical edges in $]0,n_{0}]^{d-1}\times ]0,1]$. We can thus apply the Cram\'er theorem in $\mathbb R$ (see \cite{Durrett}), which states the existence of a negative constant $c(n_{0},k_{0},\varepsilon)$ such that
$$ \lim_{m\rightarrow \infty} \frac{1}{m^{d-1}} \ln{\mathbb P\left[\frac{1}{m^{d-1}} \sum_{i_{1},...,i_{d-1}=0}^{m-1} \frac{\tau \left(\prod_{j=1}^{d-1} ]i_{j}n_{0}, (i_{j}+1)n_{0}],k_{0}\right)}{n_{0}^{d-1}} \geq \nu_{k_{0}} + \frac{3\varepsilon}{4} \right]} \, = \, c(n_{0},k_{0},\varepsilon) \, .$$
It follows that for $u=\lfloor \frac{h(n)}{2k_0} \rfloor$
\begin{align*}
\frac{1}{N^{d-1}h(N)} & \ln{\mathbb P\left[ \phi_{N^{d-1},h(N)} \geq (\nu+\varepsilon) N^{d-1} \right]} \\
& \, \leq \, \frac{u}{N^{d-1}h(N)} \ln{\mathbb P\left[ \phi_{\mathcal{B}(k_{0})} \geq (\nu_{k_{0}}+\frac{3\varepsilon}{4}) N^{d-1} \right]} \\
& \, \leq \, \frac{um^{d-1}}{N^{d-1}h(N)} \frac{1}{m^{d-1}} \ln{\mathbb P\left[\frac{1}{m^{d-1}} \sum_{i_{1},...,i_{d-1}=0}^{m-1} \frac{\tau \left(\prod_{j=1}^{d-1} ]i_{j}n_{0}, (i_{j}+1)n_{0}],k_{0}\right)}{n_{0}^{d-1}} \geq \nu_{k_{0}} + \frac{3\varepsilon}{4} \right]} \\
& \xrightarrow[m\rightarrow \infty]{} \frac{c(n_{0},k_{0},\varepsilon)}{2 k_{0}n_{0}^{d-1}} \,\, <0 \, ,\\
\end{align*}
so $\psi(\lambda)>0$. This ends the proof of theorem \ref{limite}.

\begin{rem}
The existence of a positive $\theta$ satisfying $\mathbb E [e^{\theta t}] < \infty$ is probably not a necessary condition to have the positivity of the function $\psi $ on $]\nu,+\infty[$. However, a condition on the moments of $t$ is necessary. Indeed, if the tail of the distribution of $t$ is too big, the probability to have a vertical path of edges with big capacities (bigger than $\lambda n^{d-1}$) is large, thus the probability to have $\phi_{n^{d-1},h(n)} \geq \lambda n^{d-1}$ cannot decay exponentially fast in $n^{d-1}h(n)$. 
\end{rem}


\section{Proof of Theorem \ref{pgd}}

This is an adaptation of the proof of a large deviation principle in \cite{Cerf:StFlour}.
We take $h$ such that $h(n) / \ln n \rightarrow \infty$ (we can do this
again without loss of generality because $\psi$ is independent of $h$) and we suppose that there exists a positive $\theta$ such that $\mathbb E [e^{\theta t}]$ is finite. We define
$$ \beta \, = \, \inf \{ \,v \,|\, \mathbb P[t(e)\leq v] >0 \,\} \,. $$
We remark that $\phi_{(n,...,n),h(n)}/n^{d-1}$ takes its values in
$[\beta , +\infty [$. We have to prove that
\begin{itemize}
\item for any closed subset $\mathcal{F} \subset [\beta , +\infty [ $, we have
$$ \limsup _{n\rightarrow \infty} \frac{1}{n^{d-1}h(n)} \ln \mathbb{P} \left[ \frac{\phi_{(n,...,n),h(n)}}{n^{d-1}} \in \mathcal{F} \right] \,\leq \, -\inf_{\mathcal{F}} \psi \,, $$
\item for any open subset $\mathcal{O} \subset [\beta , +\infty [ $, we have
$$\liminf _{n\rightarrow \infty} \frac{1}{n^{d-1}h(n)} \ln \mathbb{P} \left[ \frac{\phi_{(n,...,n),h(n)}}{n^{d-1}} \in \mathcal{O} \right] \,\geq \, -\inf_{\mathcal{O}} \psi \,.  $$
\end{itemize}
By definition of $\beta$, for all positive $\eta$, we have
$$ s_{\beta}(\eta) \, = \, \mathbb P [t(e) \leq \beta + \eta] \, > \, 0 \,.$$


\subsection{Upper bound}

Let $\mathcal{F}$ be a closed subset of $[\beta,+\infty[$, and $a=\inf \mathcal{F}$. Clearly
$$ \mathbb P \left[ \frac{\phi_{n^{d-1},h(n)}}{n^{d-1}} \in \mathcal{F}\right] \, \leq \, \mathbb P \left[ \frac{\phi_{n^{d-1},h(n)}}{n^{d-1}} \geq a \right] \, ,$$
so
$$ \limsup_{n\rightarrow \infty} \frac{1}{n^{d-1}h(n)} \ln \mathbb P \left[ \frac{\phi_{n^{d-1},h(n)}}{n^{d-1}} \in \mathcal{F} \right] \, \leq \, - \psi(a) \, = \, - \inf_{\mathcal{F}} \psi $$
because $\psi$ is non-decreasing on $\mathbb R^{+}$.


\subsection{Lower bound}

We shall prove the following local lower bound:
\begin{equation}
\label{local}
\forall \alpha \in [\beta,+\infty[ \,,\, \forall \varepsilon>0 \qquad \liminf_{n\rightarrow \infty} \frac{1}{n^{d-1}h(n)} \ln \mathbb P \left[ \frac{\phi_{n^{d-1},h(n)}}{n^{d-1}} \in ]\alpha-\varepsilon,\alpha+\varepsilon[ \right] \, \geq \,  - \psi(\alpha) \, .
\end{equation}
If (\ref{local}) holds, we have the desired lower bound. Indeed, if $\mathcal{O}$ is an open subset of $[\beta,+\infty[$, for every $\alpha$ in $\mathcal{O}$ there exists a positive $\varepsilon$ such that $]\alpha - \varepsilon, \alpha+\varepsilon[ \subset \mathcal{O}$, whence
\begin{align*}
\liminf_{n\rightarrow \infty} \frac{1}{n^{d-1}h(n)} \ln \mathbb P \left[\frac{\phi_{n^{d-1},h(n)}}{n^{d-1}} \in \mathcal{O} \right] & \, \geq \,  \liminf_{n\rightarrow \infty} \frac{1}{n^{d-1}h(n)} \ln \mathbb P \left[\frac{\phi_{n^{d-1},h(n)}}{n^{d-1}} \in ]\alpha-\varepsilon,\alpha+\varepsilon[ \right] \\
& \, \geq \, -\psi(\alpha) \,. \\
\end{align*}
By taking the supremum over $\alpha$ in $\mathcal{O}$, we obtain
$$ \liminf_{n\rightarrow \infty} \frac{1}{n^{d-1}h(n)} \ln \mathbb P \left[\frac{\phi_{n^{d-1},h(n)}}{n^{d-1}} \in O \right] \, \geq \, -\inf_{\mathcal{O}} \psi \,. $$
To prove (\ref{local}), we have to consider again different cases.

$\bullet \,\, \alpha \geq \nu $ : When $\psi(\alpha)=+\infty$, the result
 is obvious. For a finite $\psi(\alpha)$ we have
 $\psi(\alpha+\varepsilon)>\psi(\alpha)$ because the function $\psi$ is
 convex on $[\nu,+\infty[$, $\psi(\nu)=0$ and $\psi$ is positive on
 $]\nu,+\infty[$ so $\psi$ is increasing on $[\nu,+\infty[$ (or infinite). Now
$$ \mathbb P \left[\frac{\phi_{n^{d-1},h(n)}}{n^{d-1}} \in ]\alpha-\varepsilon,\alpha+\varepsilon[ \right] \, \geq \, \mathbb P \left[\frac{\phi_{n^{d-1},h(n)}}{n^{d-1}} \geq \alpha \right] - \mathbb P \left[ \frac{\phi_{n^{d-1},h(n)}}{n^{d-1}} \geq \alpha+\varepsilon \right] \, ,$$
so
$$ \liminf_{n\rightarrow \infty} \frac{1}{n^{d-1}h(n)} \ln \mathbb P \left[\frac{\phi_{n^{d-1},h(n)}}{n^{d-1}} \in ]\alpha-\varepsilon,\alpha+\varepsilon[ \right] \, \geq \, -\psi(\alpha) \, .$$

$\bullet \,\, \beta \leq \alpha < \nu $ :

\begin{figure}[ht!]
\centering

\begin{picture}(0,0)%
\epsfig{file=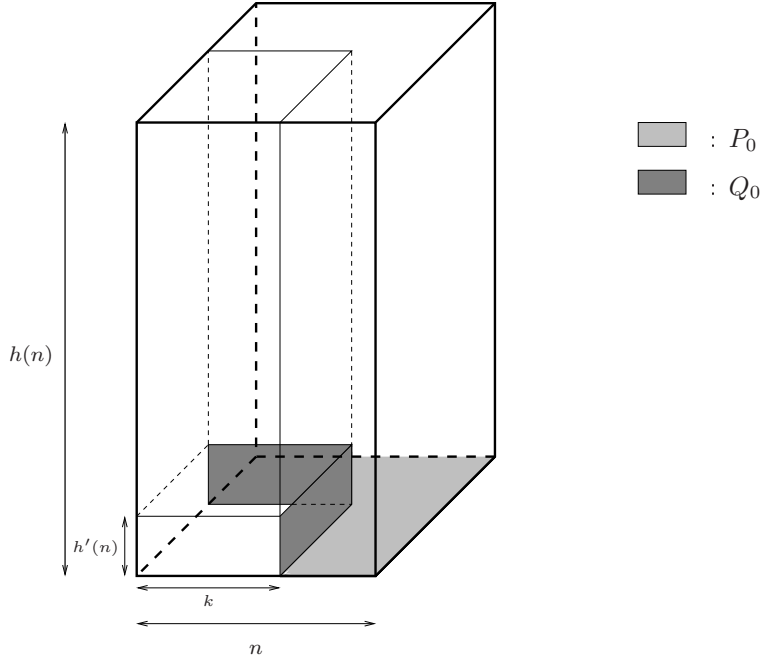}%
\end{picture}%
\setlength{\unitlength}{1973sp}%
\begingroup\makeatletter\ifx\SetFigFont\undefined%
\gdef\SetFigFont#1#2#3#4#5{%
  \reset@font\fontsize{#1}{#2pt}%
  \fontfamily{#3}\fontseries{#4}\fontshape{#5}%
  \selectfont}%
\fi\endgroup%
\begin{picture}(10034,8262)(1066,-8290)
\put(9901,-1861){\makebox(0,0)[lb]{\smash{{\SetFigFont{10}{12.0}{\rmdefault}{\mddefault}{\updefault}: $P_{0}$}}}}
\put(9901,-2461){\makebox(0,0)[lb]{\smash{{\SetFigFont{10}{12.0}{\rmdefault}{\mddefault}{\updefault}: $Q_{0}$}}}}
\put(4201,-8236){\makebox(0,0)[b]{\smash{{\SetFigFont{8}{9.6}{\rmdefault}{\mddefault}{\updefault}$n$}}}}
\put(1651,-4561){\makebox(0,0)[rb]{\smash{{\SetFigFont{8}{9.6}{\rmdefault}{\mddefault}{\updefault}$h(n)$}}}}
\put(3601,-7636){\makebox(0,0)[b]{\smash{{\SetFigFont{6}{7.2}{\rmdefault}{\mddefault}{\updefault}$k$}}}}
\put(2476,-6961){\makebox(0,0)[rb]{\smash{{\SetFigFont{6}{7.2}{\rmdefault}{\mddefault}{\updefault}$h'(n)$}}}}
\end{picture}%

\caption{Control of the flow}
\label{fig:pgd}
\end{figure}
In our cylinder $B=B((n,...,n),h(n))$ we will isolate a smaller cylinder of adequate proportions in which we will impose that the rescaled flow is around its typical value $\nu$, and we will control the amount of fluid that can circulate outside it (see figure \ref{fig:pgd}). For that purpose, we consider a function $h'$ such that
$$ h': \mathbb N \rightarrow \mathbb N \, , \,\, h'\,\leq\,h \,, \,\,
\lim_{n\rightarrow \infty}\frac{h'(n)}{n} \,=\, 0 \,\,\, and
\,\,\lim_{n\rightarrow \infty} \frac{h'(n)}{\ln n} \,=\, +\infty
\,\,\,(then\,\lim_{n\rightarrow \infty} h'(n)\, =\, +\infty) \,. $$
We define the constants
$$ v \,=\, \left( \frac{\alpha - \beta}{\nu -\beta} \right) ^{\frac{1}{d-1}} \,,\,\, k \,=\, \lfloor v n\rfloor \,,\,\, 0 \,<\, \eta \, \leq \, \frac{\varepsilon}{4} \,,$$
the set $B'$ and the corresponding event $A$
$$ B' \,=\, B((k,...,k),h(n)) \,, \,\, A\,=\, \{\, \phi_{B'} \geq (\nu-\eta)k^{d-1} \,\} \, .$$
For $i\in \mathbb N$, $0\leq i \leq (\lfloor h(n)/h'(n)\rfloor -1)$, we finally define the sets $B_{i}$, $P_{i}$ and $Q_{i}$ and the corresponding events $A_{i}$, $E_{i}$ and $F_{i}$, and the global events $E$ and $F$ as follows
$$ B_{i} \,=\, B'\, \cap\, \left( \mathbb{R}^{d-1}\times
       ]ih'(n),(i+1)h'(n)] \right) \,,$$
$$ P_{i} \,=\, (B\smallsetminus B') \cap (\mathbb Z^{d-1}\times\{\frac{1}{2} + ih'(n) \} ) \,,$$
$$ Q_{i} \,=\, \bigcup_{j=1}^{d-1} \left( [0,k]^{j-1}\times \{k+\frac{1}{2}\} \times [0,k]^{d-1-j} \times ]ih'(n),(i+1)h'(n)] \right) \,,$$
$$ A_{i}\,=\, \{\, \phi_{B_{i}} \leq (\nu+\eta)k^{d-1} \,\} \, ,$$
$$ E_{i} \,=\, \{ \, all\,\,the\,\,(n^{d-1}-k^{d-1})\,\,vertical\,\,edges\,\,e\,\,of\,\,B\smallsetminus B'\,\,that\,\,intersect\,\,P_{i}\,\,satisfy\,\,t(e)\leq\beta+\eta \,\} \,,$$
$$ E\,=\,\bigcap_{i}E_{i} \,,$$
$$ F_{i}\,=\, \{\, all\,\,the\,\,(d-1)h'(n)k^{d-2}\,\,horizontal\,\,edges\,\,e\,\,that\,\,intersect\,\,Q_{i}\,\,satisfy\,\,t(e)\leq\beta+\eta \,\} \,,$$
$$ F\,=\,\bigcap_{i}F_{i} \,.$$
Fix $n_{0}\in \mathbb N$ such that
$$ \forall n\geq n_{0} \qquad (d-1)(\beta+\eta)\frac{h'(n)}{n} \,\leq\, \frac{\varepsilon}{8} \qquad and \qquad  \left| \frac{\lfloor v n\rfloor ^{d-1}}{n^{d-1}} - v^{d-1} \right| \leq \frac{\varepsilon}{8} \,. $$
On one hand, on the event $A$, we have for $n\geq n_{0}$
\begin{align*}
 \phi_{B} & \,\geq\, n^{d-1}(\nu-\eta)\frac{\lfloor v n\rfloor ^{d-1}}{n^{d-1}} + \beta n^{d-1} \left( 1-\frac{\lfloor v n\rfloor ^{d-1}}{n^{d-1}} \right) \\
& \,\geq\, n^{d-1} \left( \nu v^{d-1} + \beta(1- v^{d-1}) - 2\frac{\varepsilon}{8}-\frac{\varepsilon}{4} \right) \\
& \,>\, n^{d-1}(\alpha - \varepsilon) \,.\\
\end{align*}
Here the term $\beta n^{d-1} ( 1- \lfloor v n\rfloor ^{d-1} / n^{d-1} )$ is
the minimal amount of fluid that crosses $B \smallsetminus B'$ from its
bottom to its top because the capacity of an edge cannot be smaller than
$\beta$, by definition of $\beta$. On the other hand, if for some $i$ in $\{ 0,..., ( \lfloor \frac{h(n)}{h'(n)} \rfloor -1 )\}$ the event $A_{i}\cap E_{i}\cap F_{i}$ occurs then we have
\begin{align*}
\forall n\geq n_{0} \qquad \phi_{B} & \,\leq\, n^{d-1} \left( (\nu+\eta) \frac{\lfloor v n\rfloor ^{d-1}}{n^{d-1}} + (\beta+\eta) \left( 1-  \frac{\lfloor v n\rfloor ^{d-1}}{n^{d-1}} + (d-1) \frac{\lfloor v n\rfloor h'(n)}{n^{d-1}} \right) \right) \\
& \,\leq\, n^{d-1} \left( \alpha + 2\frac{\varepsilon}{8} + 2\frac{\varepsilon}{4} + \frac{\varepsilon}{8} \right) \\
& \,<\, n^{d-1}(\alpha+\varepsilon) \,.\\
\end{align*}
We obtain then that
\begin{align*}
\forall n\geq n_{0} \qquad \mathbb P \left[ \frac{1}{n^{d-1}} \phi_{B} \in ]\alpha -\varepsilon,\alpha +\varepsilon[ \right] & \,\geq\, \mathbb P \left[ A\cap \left( \bigcup_{i} A_{i} \cap E_{i} \cap F_{i} \right) \right] \\
& \,\geq\, \mathbb P [E] \times \mathbb P [F] \times \mathbb P \left[ A \cap \left( \bigcup_{i} A_{i} \right) \right] \,.\\
\end{align*}
Now we know that
$$ \mathbb P [E] \,=\, s_{\beta}(\eta)^{(n^{d-1}-k^{d-1}) \lfloor \frac{h(n)}{h'(n)}\rfloor } $$
and
$$\mathbb P [F] \,=\, s_{\beta}(\eta)^{(d-1)k^{d-2}\lfloor\frac{h(n)}{h'(n)}\rfloor h'(n)} \,,$$
so
$$ \lim_{n\rightarrow\infty} \frac{1}{n^{d-1}h(n)} \ln \mathbb P [E] \,=\, \lim_{n\rightarrow\infty} \frac{1}{n^{d-1}h(n)} \ln \mathbb P [F] \,=\,0 \,.$$
Moreover we have
\begin{align*}
\mathbb P \left[ A\cap \left( \cup_{i} A_{i} \right) \right] & \,\geq\, \mathbb P [A] - \mathbb P \left[ \cap_{i} A_{i}^{c} \right] \\
& \,\geq\, \mathbb P [A] - \mathbb P [A_{0}^{c}]^{\lfloor \frac{h(n)}{h'(n)} \rfloor } \\
& \,\geq\, \mathbb P \left[ \frac{\phi_{k^{d-1},h(n)}}{k^{d-1}} \geq \nu -\eta \right] - \mathbb P \left[ \frac{\phi_{k^{d-1},h'(n)}}{k^{d-1}} \geq \nu+\eta \right]^{\lfloor \frac{h(n)}{h'(n)} \rfloor } \,,\\
\end{align*}
which leads, thanks to our previous study about $\psi$, to
$$\lim_{n\rightarrow\infty} \frac{1}{n^{d-1}h(n)} \ln \mathbb P \left[ A\cap \left( \cup_{i} A_{i} \right) \right] \,=\,0 \,.$$
We conclude that
$$\lim_{n\rightarrow\infty} \frac{1}{n^{d-1}h(n)} \ln \mathbb P \left[ \frac{\phi_{n^{d-1},h(n)}}{n^{d-1}} \in ]\alpha -\varepsilon,\alpha+\varepsilon[ \right] \,\geq\,0 \,=\, -\psi(\alpha) \,.$$
This ends the proof of the lower bound.
\\

{\bf Acknowledgments:} The author would like to warmly thank Rapha\"el Cerf
for his guidance, his help and his kindness. The author is also very
grateful to the referee for his numerous valuable comments.


\end{document}